\documentclass[11pt]{amsart} 
\usepackage{latexsym}

\usepackage{amsmath}
\usepackage{amsthm}
\usepackage{amscd}
\usepackage{amsfonts}
\usepackage[dvips]{graphicx}
\usepackage{psfrag}
\usepackage{epsfig}

\usepackage{color}
\definecolor{greencite}{rgb}{0.2,0.6,0.2} 
\definecolor{bluformula}{rgb}{0.1,0.2,0.6} 

\usepackage[colorlinks=true,linkcolor=bluformula,urlcolor=blue,citecolor=greencite]{hyperref}

\newlength{\defbaselineskip}
\newcommand{\setlinespacing}[1]
           {\setlength{\baselineskip}{#1 \defbaselineskip}}
	  
\theoremstyle{plain}
\newtheorem{theorem}{Theorem}[section]
\newtheorem{lemma}[theorem]{Lemma}
\newtheorem{proposition}[theorem]{Proposition}
\newtheorem{corollary}[theorem]{Corollary}

\newtheorem{definition}[theorem]{Definition}
\theoremstyle{definition}
\theoremstyle{remark}
\numberwithin{equation}{section}

\newcommand{\hs}{{\mathcal H}}

\newcommand{\R}{{\mathbb R}}

\newcommand{\ue}{u_{\varepsilon}}
\newcommand{\we}{w_{\varepsilon}}
\newcommand{\tue}{T{u}_{\varepsilon}}
\newcommand{\twe}{T{w}_{\varepsilon}}
\newcommand{\sie}{\sigma_{\varepsilon}}
\newcommand{\roe}{\rho_{\varepsilon}}

\newcommand{\weak}{\rightharpoonup}

\def\dys{\displaystyle}

\def\eps{\varepsilon}

\newcommand{\al}{{\alpha}}
\newcommand{\be}{{\beta}}

\newcommand{\Om}{\Omega}

\newcommand{\bvi}{BV(\Omega; \{\alpha,\beta\})}
\newcommand{\bvii}{BV(\partial\Omega; \{\alpha',\beta'\})}

\def\dys{\displaystyle}

\def\eps{\varepsilon}

\title
[$\Gamma$-convergence of some super quadratic functionals with singular weights]
{$\Gamma$-convergence of some super quadratic functionals with singular weights}
\author{Giampiero Palatucci }
\email{giampiero.palatucci@univ-cezanne.fr}
\author{Yannick Sire}
\email{sire@cmi.univ-mrs.fr}

\begin{document}
\vskip .2truecm

\subjclass[2000]{Primary  82B26, 49J45; Secondary 49Q20}

\keywords{Phase transitions, Line tension, Weighted Sobolev spaces, Nonlocal variational problems, $\Gamma$-convergence,
Functions of bounded variation}

\setlinespacing{1.1}

\begin{abstract}
\small{
We study the $\Gamma$-convergence of the following functional ($p>2$)
$$
F_{\varepsilon}(u):=\varepsilon^{p-2}\!\int_{\Omega}\!|Du|^p d(x,\partial \Omega)^{a}dx+\frac{1}{\varepsilon^{\frac{p-2}{p-1}}}\!\int_{\Omega}\!W(u) d(x,\partial \Omega)^{-\frac{a}{p-1}}dx+\frac{1}{\sqrt{\varepsilon}}\!\int_{\partial\Omega}\!V(Tu)d\mathcal{H}^2,
$$
where $\Omega$ is an open bounded set of $\mathbb{R}^3$ and $W$ and $V$ are two non-negative continuous functions vanishing at $\alpha, \beta$ and $\alpha', \beta'$, respectively.
In the previous functional, we fix $a=2-p$ and $u$ is a scalar  density function, $Tu$ denotes its trace on $\partial\Omega$, $d(x,\partial \Omega)$ stands for the distance function to the boundary $\partial\Om$.
We show that the singular limit of the energies $F_{\varepsilon}$ leads to a coupled problem of bulk and surface phase transitions.
}
\end{abstract}
\maketitle
{\small \tableofcontents}

\section{Introduction}

This paper is devoted to the $\Gamma-$convergence of the following functional $(p>2)$ 
$$
F_{\varepsilon}(u):=\varepsilon^{p-2}\!\int_{\Omega}\!|Du|^p d(x,\partial \Omega)^{a}dx+\frac{1}{\varepsilon^{\frac{p-2}{p-1}}}\!\int_{\Omega}\!W(u) d(x,\partial \Omega)^{-\frac{a}{p-1}}dx+\frac{1}{\sqrt{\varepsilon}}\!\int_{\partial\Omega}\!V(Tu)d\mathcal{H}^2,
$$
where $\Omega$ is a bounded set in $\mathbb R^3$, $V$, $W$ are two non-negative continuous functions vanishing at $\alpha,\beta$ and $\alpha',\beta'$ respectively and $a$ is a fixed number, equal to $a=2-p$; $Tu$ denotes the trace of $u$ on $\partial\Om$.

A lot of work has been devoted to the analysis of the asymptotic behavior of the functional (see for instance \cite{modica,modica87})

\begin{equation}\label{fu_momo}
\dys E_\eps(u):=\eps\int_{\Om}|{D} u|^2 dx+\frac{1}{\eps}\int_{\Om}W(u)dx.
\end{equation}

In particular, Modica proved that the previous functional $E_\varepsilon$ $\Gamma-$converges in $L^1$ to 
$$E(u)=\sigma \mathcal H^2(Su)$$
among all the admissible configurations $u \in BV(\Omega;\left \{\alpha, \beta \right \})$ with fixed volume. In the previous functional $E$, $\sigma$ is a constant depending only on the potential $W$ and $\mathcal H^2(Su)$ is the surface measure of the complement of Lebesgue points of $u$.

In \cite{alberti98}, Alberti, Bouchitt\'e and Seppecher considered the so-called two-phase model related to capillarity energy with line tension 
\begin{equation}
\dys E_\eps(u):=\eps\int_{\Om}|{D} u|^2 dx+\frac{1}{\eps}\int_{\Om}W(u)dx+\lambda_\varepsilon \int_{\partial \Omega} V(Tu) d\mathcal{H}^2.
\end{equation}
The case $\lambda_\varepsilon=\lambda$ has been considered by Modica (with $V$ being a positive continuous function), while Alberti, Bouchitt\'e and Seppecher considered a logarithmic scaling, namely $\varepsilon \log \lambda_\varepsilon \rightarrow K>0$ as $\varepsilon$ goes to $0$. 
Our approach here is to consider another penalization by perturbing with the term 
$$\int_{\Omega}\!|Du|^p d(x,\partial \Omega)^{a}dx. $$

When $a=0$, this case has been considered by one of the authors (see \cite{palatucci,palatucci07}). We consider the case when we add a weight to the gradient term, namely $d(x,\partial \Omega)$. This weight is somehow related to some non local problems involving fractional powers of the laplacian.

In the paper \cite{cafSil}, Caffarelli and Silvestre proved that one can realize any power of the fractional laplacian operator $(-\Delta)^s$ via an $s-$harmonic extension in the half-space. The fractional laplacian $(-\Delta)^s$ ($s \in (0,1)$) is a pseudo-differential operator of symbol $|\xi|^{2s}$. Caffarelli and Silvestre proved the following result:  consider the boundary Dirichlet problem
(with~$y\in \mathbb R^n$ and~$x>0$)    
\begin{equation}\label{bdyFrac} 
\left \{
\begin{matrix}
{\rm div}\, (x^a \nabla v)  =  0  \qquad 
{\mbox{ on $\mathbb R^{n+1}_+ 
:=\mathbb R^n\times(0,+\infty)$}} 
\\
v  = f, 
\qquad{\mbox{ on $\mathbb R^n\times\{0\}$,}}\end{matrix}\right . \end{equation} 
where $f$ is a given smooth compactly supported function (for instance) and $v$ is of finite energy (namely $\int_{\mathbb R^{n+1}_+} x^a |\nabla v|^2\,dx_,dy <\infty$).
Then, up to a normalizing factor,
the Dirichlet-to-Neumann operator  
$\Gamma_a:  
v|_{\partial \mathbb R^{n+1}_+} \mapsto  
-x^a v_x|_{\partial \mathbb R^{n+1}_+} $ 
is precisely $(-\Delta)^{\frac{1-a}{2}}$. 
As a consequence, one has the following corollary (see \cite{cafSil}): let $u$ be a solution of 
$$(-\Delta)^s u(y)=f(y), \ \ \ y \in \R^n
$$
and consider $P_s$ the Poisson kernel associated to the operator div$(x^{1-2s}\nabla)$. Therefore, the function $v=P_s\star_y u$ is a solution of the following problem 
\begin{equation}\label{bdyFrac2} 
\left \{
\begin{matrix}
{\rm div}\, (x^a \nabla v)  =  0  \qquad 
{\mbox{ on $\mathbb R^{n+1}_+ 
:=\mathbb R^n\times(0,+\infty)$}} 
\\
v  = u, 
\qquad{\mbox{ on $\mathbb R^n\times\{0\}$,}}
\\
-x^av_x=f\qquad{\mbox{ on $\mathbb R^n\times\{0\}$.}}
\end{matrix}\right . \end{equation} 
Note that  the condition~$ \frac{1-a}{2}=s\in(0,1)$  
reduces to~$a \in(-1,1)$. The weight $x^a$ is a particular weight since it belongs to Muckenhoupt $A_2$ classes (see \cite{muck}). Indeed, since $a \in (-1,1)$, the weight $x^a$ (as its inverse) is locally integrable.  

A quick look at the weight $x^a$ shows that it is just the distance of a point $(x,y) \in \mathbb R^{n+1}_+$ to the boundary of the domain, namely $\partial \mathbb R^{n+1}_+ =\mathbb R^n.$ Therefore, a natural generalization in bounded domains consists in taking as the weight the distance to the boundary $d(x,\partial \Omega)^a$. In this case, there are no results available to describe what is precisely the boundary operator. However, one can expect that such a weight produces some new geometrical effects. 

In the present work, we concentrate on a quasi-linear functional $F_{\varepsilon}$, i.e. $p>2$. The case $p=2$ has been considered in \cite{gonzalez}. In this case, the main point consists in replacing the penalizing term of the functional by its Sobolev trace norm. To be able to do such a trick, which goes back to \cite{alberti}, one has to consider the optimal Sobolev embedding, i.e. to use the optimal constant in the Sobolev inequality. Using Caffarelli-Silvestre extension technique, Gonzalez computed explicitely the constant of this embedding.

The case $p>2$ involves more technicalities due to the quasi-linear feature of the perturbation. In particular, we do not know how to replace the penalizing term by a Sobolev trace norm. Another difficulty comes from the scaling property. Indeed, 
the super-quadratic case enjoys a natural scaling which forces the parameter $a$ in the functional to be exactly $a=2-p$. As a consequence, as soon as $p \geq 3$, the weight $d(x,\partial \Omega)^a$ is no longer locally integrable. In this case, Nekvinda (see \cite{nekvinda}) proved that functions of the weighted Sobolev space $W^{1,p}(\Omega, d(x,\partial \Omega)^a)$ have no trace on $\partial \Omega$. Therefore, one has to use new techniques to deal with this case. As a consequence of this, we will be constrained to the range $p \in (2,3).$ 

To simplify notations, we will denote $h(x)=d(x,\partial \Omega)$ and then consider the following functional 
    
\begin{equation}\label{funzionale}
\dys F_{\varepsilon}(u):=\varepsilon^{p-2}\!\int_{\Omega}\!|Du|^p h^{2-p}dx+\frac{1}{\varepsilon^{\frac{p-2}{p-1}}}\!\int_{\Omega}\!W(u) h^{\frac{p-2}{p-1}}dx+\frac{1}{\sqrt{\varepsilon}}\!\int_{\partial\Omega}\!V(Tu)d\mathcal{H}^2. 
\end{equation}
\bigskip

\medskip

At this point, some remarks on the scaling have to be noticed. Choosing $\eps^{\frac{p-2}{p-1}}$ to denote the length of the bulk transition, by standard scaling analysis, the power $\eps^{p-2}$ follows naturally in the perturbation term. The election of the square root of $1/\eps$ in the boundary term is justified by the scaling property of the functional $F_\eps$ (see Section \ref{ch_2d}).

\section{Description of the results}\label{sec_descri}

We first fix notations, recalling also some standard mathematical results used throughout the paper. Then, we analyze the asymptotic behavior of the functional $F_\eps$ defined in (\ref{funzionale}) stating the related main convergence result.

\medskip
\subsection{Notation}

In this work, we consider different domains $A$ in dimensions $n=1, 2, 3$; more precisely, $A$ will always be a bounded open set of $\R^n$. 
We denote by $\partial A$ the boundary of $A$ relative to the ambient space; $\partial A$ is always assumed to be Lipschitz regular.
 Unless otherwise stated, $A$ is  endowed with the corresponding $n$-dimensional Hausdorff measure, $\hs^n$ (see \cite{evans}, Chapter 2). We write $\dys \int_{A}\!fdx$ instead of $\dys \int_{A}\! fd\hs^n$. 

The {\em essential boundary} of $A$ is the set of all points where $A$ has neither 0 nor 1 density and where the density does not exist. Since the essential boundary agrees with the topological boundary when the latter is Lipschitz regular, we also denote the essential boundary by $\partial A$.

\bigskip

For every $u \in L^1_{\text{loc}}(A)$, we denote by $Du$ the derivative of $u$ in the sense of distributions. As usual, for every $p\geq 1$, $W^{1,p}(A)$ is the Sobolev space of all $u\in L^p(A)$ such that $Du \in L^p(A)$.
Given a weight $w:A\to[0,\infty)$, and $p\geq 1$, we consider the weighted Sobolev space  $W^{1,p}(A, w)$ the space of all functions $u$ with norm
$$
\dys \|u\|^p_{W^{1,p}(A, w)} := \int_{A}|u|^p w dx +\int_{A} |Du|^p w dx.
$$

$BV(A)$ is the space of all $u\in L^1(A)$ with bounded variation; i.e., such that $Du$ is a bounded Borel measure on $A$. We denote by $Su$ the {\em jump set}; i.e., the complement of the set of Lebesgue points of $u$.

\bigskip

For every $s\in (0,1)$ and every $p\geq 1$, $W^{s,p}(A)$ is the space of all $u\in L^p(A)$ such that the fractional semi-norm $\dys \int_A\!\!\int_{A}\frac{|u(x)-u(x')|^p}{|x-x'|^{sp+n}}dxdx'$ is finite.

We denote by $T$ the trace operator which maps $BV(A)$ onto $L^1(\partial A)$ and $W^{1,p}(A, w)$ onto $W^{2-3/p,p}(\partial A)$, for a suitable weight $w$ (see \cite[Theorem 2.8]{nekvinda}). In particular for $p\in (2,3)$, there exists a constant $S_p$ such that
$$
\dys \|Tu\|_{W^{2-3/p,p}(\partial\Om)} \leq S_p\|u\|_{W^{1,p}(\Om, d^{2-p}\!(x,\partial\Om))} \ \ \ \text{(see \cite[Theorem 2.11]{nekvinda})}.
$$
For details and results about the theory of $BV$ functions and Sobolev spaces we refer to \cite{evans}, \cite{ambrosio} and \cite{adams}.

\medskip

\subsection{The $\Gamma$-convergence result}\label{sec_profilep}

Let $\Om$ be a bounded open subset of $\R^3$ with smooth boundary; let $W$ and $V$ be non-negative continuous functions on $\R$ with growth at least linear at infinity and vanishing respectively only in the ``double well'' $\{\alpha, \beta\}$, with $\alpha<\beta$, and $\{\alpha', \beta'\}$, with $\alpha'<\beta'$.
Assume that the potential $V$ is convex near its wells.

Let $p \in (2,3)$ be a real number. For every $\eps>0$ we consider the functional $F_\eps$ defined in ${W^{1, p}({\Omega}, h^{2-p})}$, given by
\begin{equation}\label{funzionale1}
F_\eps(u):=\eps^{p-2}\int_\Om|{D} u|^p h^{2-p} dx + \frac{1}{\eps^{\frac{p-2}{p-1}}}\int_{\Om}W(u) h^{\frac{p-2}{p-1}}dx+\frac{1}{\sqrt{\eps}}\int_{\partial\Om}V(Tu)d\hs^2.
\end{equation}

\smallskip

We analyze the asymptotic behavior of the functional $F_\eps$ in terms of $\Gamma$-convergence. Let $(\ue)$ be an equi-bounded sequence for $F_\eps$; i.e., there exists a constant $C$ such that $\dys F(\ue)\leq C$.
We observe that the term $\dys \frac{1}{\eps^{\frac{p-2}{p-1}}}\int_{\Om}W(u_\eps) h^{\frac{p-2}{p-1}}dx$ forces $\ue$ to take values close to $\alpha$ and $\beta$, while the term $\dys \eps^{p-2}\int_\Om|{D} u_\eps|^p h^{2-p}dx$ penalizes the oscillations of $\ue$. We will see that when $\eps$ tends to $0$, the sequence $(\ue)$ converges (up to a subsequence) to a function $u$, belonging to $BV(\Om)$, which takes only the values $\alpha$ and $\beta$. Moreover each $\ue$ has a transition from the value $\alpha$ to the value $\beta$ in a thin layer close to the surface $Su$, which separates the bulk phases $\{u=\alpha\}$ and $\{u=\beta\}$.
Similarly, the boundary term of $F_\eps$ forces the traces $\tue$ to take values close to $\alpha'$ and $\beta'$, and the oscillations of the traces $\tue$ are again penalized by the integral $\dys \eps^{p-2}\int_\Om|{D} u_\eps|^p h^{2-p}dx$. Then, we expect that the sequence $(\tue)$ converges to a function $v$ in $BV(\partial\Om)$ which takes only the values $\alpha'$ and $\beta'$, and that a concentration of energy occurs along the line $Sv$, which separates the boundary phases $\{v=\alpha'\}$ and $\{v=\beta'\}$.
\smallskip

In view of possible ``dissociation of the contact line and the dividing line'' (see \cite[Example 5.2]{alberti98}), we recall that $Tu$ may differ from $v$. Since the total energy $F_\eps(\ue)$ is partly concentrated in a thin layer close to $Su$ (where $\ue$ has a transition from $\alpha$ to $\beta$), partly in a thin layer close to the boundary (where $\ue$ has a transition from $Tu$ to $v$), and partly in the vicinity of $Sv$ (where $\tue$ has a transition from $\alpha'$ to $\beta'$), we expect that the limit energy is the sum of a surface energy concentrated on $Su$, a boundary energy on $\partial\Om$ (with density depending on the gap between $Tu$ and $v$), and a line energy concentrated along $Sv$.

The asymptotic behavior of the functional $F_\eps$ is described by a functional $\Phi$ which depends on the two functions $u$ and $v$.
Let $\mathcal{W}$ be an antiderivative of $W^{(p-1)/p}$.
 For every $(u,v) \in \bvi\times\bvii$, we will prove that
\begin{equation}\label{fu_phi}
\dys \Phi(u,v):=\sigma_p\hs^2(Su)+c_p\int_{\partial\Om}|\mathcal{W}(Tu)-\mathcal{W}(v)| d\hs^2+\gamma_p\hs^1(Sv),
\end{equation}
where as usual the jump sets $Su$ and $Sv$ are the complement of the set of Lebesgue points of $u$ and $v$, respectively; $\dys c_p$ and $\sigma_p$ are the constants defined by
\begin{eqnarray}\label{def_cp}
\dys c_p:=\frac{p}{(p-1)^{(p-1)/p}}, \,\,\, 
\sigma_p:= c_p|\mathcal{W}(\beta)-\mathcal{W}(\alpha)|; 
\end{eqnarray}
The constant $\gamma_p$ is given by the optimal profile problem
\begin{eqnarray}\label{def_gamma-z}
\gamma_{p}:=\inf\left\{\int_{\R^2_+}|{D} u|^p x_2^{2-p} dx+\int_{\R}V(Tu)d\hs^1 : u \in L^1_{\text{loc}}(\R^2_+) : \hskip 2cm \right.  \\
\left. \lim_{t\to-\infty}\!Tu(t)\!=\!\alpha',  \ \lim_{t\to+\infty}\!Tu(t)\!=\!\beta'\right \}.
\end{eqnarray}
Note that in the definition (\ref{def_gamma-z}) we utilize the variables $x=(x_1, x_2)\in \mathbb R \times \mathbb R^+$ to denote any point of $\R^2_+$, so that we have always $h^{2-p}=x_2^{2-p}$. \medskip

The main convergence result is precisely stated in the following theorem.
\begin{theorem}\label{c4_teoremap}
Assume $p \in (2,3)$. Let $F_{\eps}:W^{1,p}(\Om, h^{2-p})\to\R$ and $\Phi:\bvi\times\bvii \to\R$ defined by {\rm (\ref{funzionale1})} and {\rm (\ref{fu_phi})}.

Then

\begin{itemize}
\item[{\rm (i)}]{\textsc{[Compactness]} If $(\ue) \subset {W^{1, p}({\Omega}, h^{2-p})}$ is a sequence such that $F_{\eps}(\ue)$ is bounded, then $(\ue,\tue)$ is pre-compact in $L^{1}(\Om)\times L^1(\partial\Om)$ and every cluster point belongs to $\bvi\times\bvii$.}
\item[{\rm (ii)}]{\textsc{[Lower Bound Inequality]} For every $(u,v) \in\bvi\times\bvii$ and every sequence $(\ue)\subset{W^{1, p}({\Omega}, h^{2-p})}$ such that $\ue \to u$ in $L^{1}(\Om)$ and $\tue \to v$ in $L^{1}(\partial\Om)$,
\begin{equation*}
\liminf_{\eps\to0} F_{\eps}(\ue)\geq \Phi(u,v).
\end{equation*}}
\item[{\rm (iii)}]{\textsc{[Upper Bound Inequality]} For every $(u, v) \in\bvi\times\bvii$ there exists a sequence $(\ue)\subset{W^{1, p}({\Omega}, h^{2-p})}$ such that $\ue \to u$ in $L^{1}(\Om)$, $\tue \to v$ in $L^{1}(\partial\Om)$ and
\begin{equation*}
\limsup_{\eps\to0} F_{\eps}(\ue)\leq \Phi(u,v).
\end{equation*}}
\end{itemize}
\end{theorem}

We can easily rewrite this theorem in term of $\Gamma$-convergence. To this aim, we extend each $F_{\eps}$ to $+\infty$ on $L^{1}(\Om)\setminus {W^{1, p}({\Omega}, h^{2-p})}$ and, from Theorem \ref{c4_teoremap}, we deduce that 

\begin{corollary}\label{c4_corollario}
$F_\eps \ \Gamma$-converges on $L^1(\Om)$ to $F$, given by
\begin{equation*}
\dys F(u):=\begin{cases} \inf\left\{\Phi(u,v): v\in \bvii\right\} & \text{if} \ u \in \bvi, \\
+\infty & \text{elsewhere in} \ L^1(\Om).
\end{cases}
\end{equation*}
\end{corollary}

\bigskip

\section{Strategy of the proof and some convergence results}\label{sec_strategy}

 The proof of Theorem \ref{c4_teoremap} requires several steps in which we have to analyze different effects. Then, we can deduce the terms of the limit energy $\Phi$, localizing three effects: the bulk effect, the wall effect and the boundary effect.

\subsection{The bulk effect}
In the bulk term, the limit energy can be evaluated like in \cite{gonzalez}. This requires to generalize the Modica-Mortola results on the functional (\ref{fu_momo}) (see \cite{momo}) to a functional with super-quadratic growth in the perturbation term involving the singular weight $h^{2-p}$.
\medskip

For every open set $A\subset\R^3$, $p \in (2,3)$ and every real function $u\in W^{1,p}(A, h^{2-p})$, we consider the functional
\begin{equation}\label{fu_momow}
\dys G_\eps(u, A):=\eps^{p-2}\int_{A}|{D} u|^p h^{2-p} dx+\frac{1}{\eps^{\frac{p-2}{p-1}}}\int_{A}W(u) h^{\frac{p-2}{p-1}}dx.
\end{equation}
Since there is no interaction with the boundary of $A$ and the weight is regular in the interior, the asymptotic behavior of the functional $G_{\eps}$ will be very similar to the one of (\ref{fu_momo}).
\begin{theorem}\label{teo_modicap}
For every domain $A\subset\Om$ the following statements hold.
\begin{itemize}
\item[{\rm (i)}]{If $(\ue)\subset W^{1,p}(A, h^{2-p})$ is a sequence with uniformly bounded energies $G_\eps(\ue, A)$. Then $(\ue)$ is pre-compact in $L^1(A)$ and every cluster point belongs to $BV(A; \{\al,\be\})$.}
\item[{\rm(ii)}]{For every $u\in BV(A; \{\al,\be\})$ and every sequence $(\ue)\subset W^{1,p}(A, h^{2-p})$ such that $\ue\to u$ in $L^1(A)$,
$$
\liminf_{\eps\to 0} G_\eps(\ue, A) \geq \sigma_p\hs^2(Su\cap A),
$$}
\item[{\rm (iii)}]{For every $u\in BV(A; \{\al,\be\})$ there exists a sequence $(\ue)\subset W^{1,p}(A)$ such that $\ue\to u$ in $L^1(A)$ and
$$
\limsup_{\eps\to 0} G_\eps(\ue, A) \leq \sigma_p\hs^2(Su\cap A).
$$
}
\end{itemize}
\end{theorem}

\noindent
{\bf Proof.} The proof is close to the one of Gonzalez in \cite[Proposition 3.1]{gonzalez} and Modica-Mortola's one. Here we provide a sketch and the needed modifications due to the different growth power in the singular perturbation.
\smallskip

Using the following Young's inequality, $X, Y \geq 0$,
\begin{equation}\label{eq_young}
XY \leq \dys \frac{X^p}{p}+\frac{Y^q}{q} , \ \ \ \bigl(q : {1}/{p}+{1}/{q}=1\bigl),
\end{equation}
with
$$
\dys X= |Du|h^{\frac{2-p}{p}}p^{\frac{1}{p}}\eps^{-\frac{2-p}{p}} \ \ \text{and} \ \ Y=W(u)^{\frac{1}{q}}h^{-\frac{2-p}{(p-1)q}}q^{\frac{1}{q}}\eps^{-\frac{p-2}{(p-1)q}},
$$
we obtain
\begin{equation}\label{lb_bulk}
\dys G_{\eps}(u,A)\geq c_p\int_A W^{\frac{p-1}{p}}|Du|dx = c_p\int_A |D(\mathcal{W}(u)|dx,
\end{equation}
where $\mathcal{W}$ is a primitive of $W^{(p-1)/p}$ and $c_p$ is defined by (\ref{def_cp}). This gives the compactness result (i) and the lower bound inequality (ii), using standard arguments.
\medskip

Consider a function $u$ in $BV(A; \{\al, \be\})$. To construct the recovery sequence $u_\eps$ of the upper bound inequality (iii), we need to take care of the weight $h^{2-p}$.
\smallskip

First, without loss of generality, we may assume that the singular set $Su$ of $u$ is a Lipschitz surface in $A$ (\cite[Theorem 1.24]{giusti}). For every $x$ in $A$, let us define the signed distance from $Su$ as
$$
\dys d'(x):=\begin{cases} \text{dist}(x, Su) &\text{if} \ x \in \{ u=\beta\}, \\
-\text{dist}(x, Su) &\text{if} \ x \in \{ u=\alpha\}.
\end{cases}
$$
We may consider smooth coordinates $(d'(x),\eta)$ in $A$ such that $\eta$ parametrizes $Su$.

Now, we choose $\theta \in W^{1,1}_{\text{\rm loc}}(\R)$ satisfying
\begin{equation}\label{eq_profile}
\dys 
\begin{cases}
\theta'=  \frac{1}{(p(p-1))^{1/p}}W^{1/{p}}(\theta) \ \  \text{a.e.} \\
\\
\theta(-\infty)=\al, \ \theta(+\infty)=\be,
\end{cases}
\end{equation}
where the valus $\theta(\pm\infty)$ are understood as the existence of the corresponding limits.
We remark  that this real function $\theta$ is just the optimum profile for the case $a=0$.

Consider the function $\phi:A\to\R$ defined by
\begin{equation}\label{fu_phisigma}
\dys \phi(t,\eta)\equiv \phi_\eta(t):=\theta\left(\frac{t}{h^{\frac{2-p}{p-1}}(0,\eta)}\right).
\end{equation}

Finally, we are in position to  construct the recovery sequence $(\ue)$. For every $\eps>0$, let $t=d'(x)/\eps^{\frac{p-2}{p-1}}$ and
\begin{equation*}
\dys \ue(x):=\phi_\eta\left(\frac{d'(x)}{\eps^{\frac{p-2}{p-1}}}\right) \ \ \forall x\in A.
\end{equation*}
Using the fact that for every $\delta \in (0,1)$ there exists $c(\delta)\to\infty$ when $\delta\to 0$ such that
$$
(X+Y)^p\leq (1+\delta)X^p+c(\delta)Y^p,
$$
by definition of $\ue$ we have
\begin{eqnarray}\label{fu_diff}
\dys |D\ue|^p(x) \! & = & \! \left| \frac{\partial\phi}{\partial t}(t(x),\eta)Dt(x)+\frac{\partial\phi}{\partial\eta}(t(x),\eta)\right|^p \nonumber \\
\\
& \leq & \! (1+\delta)\frac{(\phi'_\eta(t))^p}{\eps^{\frac{p-2}{p-1}p}}+c(\delta)\mathcal{R}(\eta,t) \ \ \forall x \in A, \nonumber
\end{eqnarray}
where we denoted by $\dys \mathcal{R}(\eta,t):=\left|\frac{\partial\phi}{\partial\eta}(t,\eta)\right|^p$.

\smallskip

Thus we can estimate the energy of the function $\ue$, using the CoArea Formula. For every $\delta\in (0,1)$ we have
\begin{eqnarray*}
\dys
G_\eps(\ue, A) \! & = & \! \eps^{p-2}\int_{A}|{D} u|^p h^{2-p} dx+\frac{1}{\eps^{\frac{p-2}{p-1}}}\int_{A}W(u) h^{\frac{p-2}{p-1}}dx \\
\\
& \leq  & \! (1+\delta)\frac{1}{\eps^{\frac{p-2}{p-1}}}\int_A \left[ (\phi'_{\eta}(t))^p h^{2-p} + W(\phi_{\eta}(t))h^{\frac{p-2}{p-1}} + c(\delta)\mathcal{R}(\eta,t)\eps^{\frac{(p-2)p}{p-1}}\right]dx \\
\\
& = & \! (1+\delta)\int_{-\infty}^{+\infty}\!\int_{\Sigma_{\eps}}\left[(\phi'_{\eta}(t))^p h^{2-p}+W(\phi_{\eta}(t))h^{\frac{p-2}{p-1}} + c(\delta)\mathcal{R}(\eta,t)\eps^{\frac{(p-2)p}{p-1}}\right]d\eta dt,
\end{eqnarray*}
with the level set $\dys \Sigma_\eps:=\left\{ x\in A : d(x, Su)=\eps^{\frac{p-2}{p-1}}t\right\}$ that converges to $Su\cap A$ when $\eps\to 0$. Moreover, when $\eps$ goes to 0,  $c(\delta)\mathcal{R}(\eta,t)\eps^{\frac{(p-2)p}{p-1}}$ converges to 0 and, if $x$ is written in the coordinates $(d'(x),\eta)$, then $h(t,\eta)$ converges to $\text{dist}((0,\eta),\partial A)\equiv h(0,\eta)$.
Hence, for every $\delta\in (0,1)$, taking the limit as $\eps$ goes to 0, we have
\begin{equation}\label{eq_quasilimsup}
\dys
\limsup_{\eps\to0} G_\eps(\ue,A) \leq (1+\delta)\int_{-\infty}^{\infty}\!\int_{Su\cap A}\left[(\phi'_{\eta}(t))^p h^{2-p}(0,\eta) + W(\phi_\eta(t))h^{\frac{p-2}{p-1}}(0,\eta)\right]d\eta dt.
\end{equation}\smallskip
Using the definitions of $\theta$ and $\phi_\eta$ given by (\ref{eq_profile}) and (\ref{fu_phisigma}), it follows that
$$
\dys \phi'_\eta(t)h^{\frac{2-p}{p}}(0,\eta)p^{\frac{1}{p}}=\left(W(\phi_\eta(t))^{\frac{1}{q}}h^{-\frac{2-p}{(p-1)q}}q^{\frac{1}{q}}\right)^\frac{1}{p-1}.
$$
So, when we apply the inequality (\ref{eq_young}), like in (\ref{lb_bulk}), we also have an equality and the (\ref{eq_quasilimsup}) becomes
\begin{eqnarray*}
\dys 
\limsup_{\eps\to0} G_\eps(\ue,A) \! & \leq & \! (1+\delta)\int_{Su\cap A}\!\int_{-\infty}^{+\infty}\! c_pW^{\frac{p-1}{p}}(\phi_{\eta}(t))\phi'_{\eta}(t) dt d\eta \\
\\
& = & \! (1+\delta)\int_{Su\cap A}\!\int_{\al}^{\be}c_p W^{\frac{p-1}{p}}(r)dr d\eta \\
\\
& = & (1+\delta)\sigma_p\hs^2(Su\cap A) \ \ \ \forall \delta\in (0,1).
\end{eqnarray*}
This concludes the proof. \hfill $\Box$

\medskip

\subsection{The wall effect}

The second term of $\Phi$ can be obtained thanks to the following lemma.
\begin{proposition}\label{pro_modica2p} 
For every domain $A\subset\Om$ with boundary piecewise of class $C^1$ and for every $A'\subset \partial A$ with Lipschitz boundary, the following statements hold.
\begin{itemize}
\item[{\rm (i)}]{For every $(u, v) \in BV(A; \{\al,\be\})\times BV(A'; \{\al',\be'\})$ and every sequence $(\ue) \subset W^{1,p}(A, h^{2-p})$ such that $\ue \to u$ in $L^1(A)$ and $\tue \to v$ in $L^1(A')$,
$$
\dys \liminf_{\eps\to 0}G_\eps(\ue, A)\geq c_p\int_{A'}|\mathcal{W}(Tu)-\mathcal{W}(v)|d\hs^2.
$$}
\item[{\rm (ii)}]{Let a function $v$, constant on $A'$, and a function $u$, constant on $A$, such that $u\equiv\alpha$ or $u\equiv\beta$, be given. Then there exists a sequence $(\ue)$ such that $\tue=v$ on $A'$, $\ue$ converges uniformly to $u$ on every set with positive distance from $A'$ and
$$
\dys \limsup_{\eps\to 0}G_\eps(\ue, A)\leq c_p\int_{A'}|\mathcal{W}(Tu)-\mathcal{W}(v)|d\hs^2.
$$
}
\end{itemize}
Moreover, the function $\ue$ may be required to be $C_r$-Lischitz continuous in $A_r:=\bigl\{x\in A : d(x, \partial A) \leq r \bigl\}$.
\end{proposition}
\smallskip
\noindent
\\ {\bf Proof.} The proof of (i) is essentially contained in \cite[Proposition 1.2 and Proposition 1.4]{modica87}, where Modica study a Cahn-Hilliard functional with quadratic growth in the singular perturbation term and with a boundary contribution (see also \cite[Proposition 4.3]{palatucci} for details of the super-quadratic version).
While, the proof of (ii) is very similar to \cite[Proposition 4.3]{alberti98} and it can be recovered using the modifications introduced in the proof of Theorem \ref{teo_modicap}-(iii).
See also \cite[Proposition 3.1]{gonzalez} for the computation of the Lipschitz constant of $\ue$.
\medskip

\subsection{The boundary effect}
This is a delicate step, that requires a deeper analysis. The main strategy is the one used by Alberti, Bouchitt\'e and Seppecher in \cite{alberti98} with the needed modifications introduced by one of the author in \cite{palatucci07} for functionals with super-quadratic growth in the singular perturbation term.
We reduce to the case in which the boundary is flat; hence we study the behavior of the energy in the three-dimensional half ball; then we reduce the problem to one dimension via a slicing argument.

Thus, the main problem becomes the analysis of the asymptotic behavior of the following two-dimensional functional
\begin{equation}\label{fu_h2}
\dys H_\eps(u):=\eps^{p-2}\int_{D_1}|{D} u|^p x_2^{2-p} dx+\frac{1}{\sqrt{\eps}}\int_{E_1}V(Tu)d\hs^1, \ \ \ \forall u \in W^{1,p}(D_1, h^{2-p}),
\end{equation}
where $D_1$ and $E_1$ are defined by
\begin{eqnarray}\label{def_drer}
&\dys D_r:=\left\{(x_1,x_2)\in\R^2: x_1^2+x_2^2<r^2, x_2>0\right\}, \nonumber \\
\\
&E_r:=\left\{(x_1,x_2)\in\R^2: x_1^2+x_2^2<r, x_2=0\right\}\equiv(-r,r).\nonumber
\end{eqnarray}

\medskip

Note that for the quadratic case the two-dimensional Dirichlet weighted energy can be replaced on the half-disk $D_r$ by the $H^{1/2}$ intrinsic norm on the ``diameter'' $E_r$. This is possible thanks to the existence of an optimal constant for the trace inequality involving the weighted $L^2$-norm of the gradient of a function defined on a two-dimensional domain and the $H^{1/2}$-norm of its trace on a line (see \cite[Proposition 4]{gonzalez}). Hence, the analysis of the line tension effect is reduced to the one of the following one-dimensional perturbation problem involving a non-local term:
$$
\dys E_\eps(v)=\eps^{1-a}\!\int\!\!\int_{I\times I}\!\frac{|v(t)-v(t')|^2}{|t-t'|^{1+2s}}dt'dt+\frac{1}{\eps^{\frac{1-a}{-a}}}\!\int_{I}\! V(v)dt, \ \ \ (I \ \text{open interval of} \ \R; \ s=(1-a)/{2}),
$$
that was essentially studied by Garroni and one of the author in \cite{garroni06}.

On the contrary, we have to study the asymptotic analysis of $H_\eps$, that will be the subject of Section \ref{ch_2d}.

\bigskip

We conclude this section stating some properties of the functional $F_\eps$.
	
	\medskip
	
\subsection{Some remark about the structure of $F_\eps$}\label{sec_structure}

The methods used in the proof of the main results of this paper strongly requires the ``localization'' of the functional $F_\eps$; i.e., looking at $F_\eps$ as a function of sets. By fixing $u$ we will be able to characterize the various effects of the problem.
In this sense, for every open set $A\subset\R^3$, every set $A'\subseteq\partial A$ and every function $u\in W^{1,p}(A, h^{2-p})$, we will denote
\begin{equation}\label{fu_loc}
\dys F_\eps(u, A, A'):=\eps^{p-2}\int_A|{D} u|^p h^{2-p} dx + \frac{1}{\eps^{\frac{p-2}{p-1}}}\int_{A}W(u) h^{\frac{p-2}{p-1}} dx+\frac{1}{\sqrt{\eps}}\int_{A'}V(Tu) d\hs^2.
\end{equation}
Clearly, $F_\eps(u)=F_\eps(u,\Om,\partial\Om)$ for every $u\in{W^{1, p}({\Omega}, h^{2-p})}$.
\smallskip

Let us observe that, thanks to the growth hypothesis on the potentials $W$ and $V$, we may assume that there exists a constant $m$ such that:
\begin{eqnarray}\label{ipotesi}
&\ -m\leq \alpha, \alpha', \beta, \beta' \leq m,  \nonumber\\
&\dys W(t)\geq W(m) \ \text{and} \ V(t)\geq V(m) \ \text{for} \ t\geq m,\\
&\dys W(t)\geq W(-m) \ \text{and} \ V(t)\geq V(-m) \ \text{for} \ t\leq -m.\nonumber
\end{eqnarray}
In particular, assumption (\ref{ipotesi}) will allow us to use the truncation argument given by the following Lemma.
\begin{lemma}\label{lem_tronca}
Let a domain $A\subset\R^3$, a set $A'\subseteq\partial A$, and a sequence $(\ue)\subset W^{1,p}(A, h^{2-p})$ with uniformly bounded energies $F_{\eps}(\ue, A, A')$ be given.
 \\ If we set ${\bar{u}_{\eps}}(x):=\max\{\min\{\ue(x), m\},-m\}$, then
\begin{itemize}
\item[{\rm (i)}]{$\dys F_\eps({\bar{u}_{\eps}}, A, A')\leq F_\eps(\ue, A, A')$,}
\item[{\rm (ii)}]{$\|{\bar{u}_{\eps}}-\ue\|_{L^1(A)}$ and $\|T{\bar{u}_{\eps}}-\tue\|_{L^1(A')}$ vanish as $\eps\to0$.}
\end{itemize}
\end{lemma}
\medskip
\noindent
\\ {\bf Proof.} 
The inequality $\dys F_\eps({\bar{u}_{\eps}}, A, A')\leq F_\eps(\ue, A, A')$ follows immediately from (\ref{ipotesi}). Statement (ii) follows from the fact that both $W$ and $V$ have growth at least linear at infinity and the integrals $\dys \int W(\ue) h^{\frac{p-2}{p-1}}dx$ and $\dys \int V(\tue)d\hs^2$ vanish as $\eps$ goes to 0.
This is a standard argument; see, for instance, \cite[Lemma 4.4]{palatucci}. \hfill $\Box$

\bigskip

    \section{Recovering the ``contribution of the wall'': the flat case}\label{ch_2d}
    
\ \ \ \ We will obtain ``the contribution of the wall''  to the limit energy $\Phi$, defined by (\ref{fu_phi}), namely $\gamma_p\hs^1(Sv)$, by estimating the asymptotic behavior of the functional
$$
\dys F_\eps(u,B\cap\Om, B\cap\partial\Om)\!=\!\eps^{p-2}\!\int_{B\cap\Om}\!|Du|^ph^{2-p}dx+\frac{1}{\eps^{\frac{p-2}{p-1}}}\int_{B\cap\Om}\!W(u)h^{\frac{p-2}{p-1}}dx+\frac{1}{\sqrt{\eps}}\!\int_{B\cap\partial\Om}\!V(Tu)d\hs^2,
$$
when $B$ is a small ball centered on $\partial\Om$ and $B\cap\partial\Om$ is a flat disk. We will follow the idea of Alberti, Bouchitt\'e and Seppecher in \cite{alberti98}, using a suitable slicing argument; the flatness assumption on $B\cap\partial\Om$ can be dropped when $B$ is sufficiently small. Hence, we need to prove a compactness result and a lower bound inequality for the following two-dimensional functional
\begin{equation}\label{def_fu2d}
\dys H_\eps(u):=\eps^{p-2}\int_{D_1}|{D} u|^p h^{2-p} dx+\frac{1}{\sqrt{\eps}}\int_{E_1}V(Tu)d\hs^1, \ \ \ \forall u \in W^{1,p}(D_1, h^{2-p}; [-m,m]),
\end{equation}
where $E_1$ and $D_1$ are defined by (\ref{def_drer}). We recall that we will always study $H_\eps$ like a reduction of $F_\eps$. Hence there will be some hypotheses inherited by this reduction. In particular, the hypothesis $u\in[-m,m]$ in (\ref{def_fu2d}) is justified by Lemma \ref{lem_tronca}.

Let us introduce the ``localization'' of the functional $H_\eps$. For every open set $A\subset\R^2$, every set $A'\subset\partial A$ and every function $u\in W^{1,p}(A, h^{2-p})$, we will denote
\begin{equation}\label{fu_disco}
\dys H_\eps(u, A, A'):=\eps^{p-2}\int_{A}|{D} u|^p h^{2-p}  dx+\frac{1}{\sqrt{\eps}}\int_{A'}V(Tu)d\hs^1.
\end{equation}

Let $A=D_1$ be the half disk defined in (\ref{def_drer}) and denote by
\begin{equation}\label{def_d10}
\dys D_1^0:=\{x=(x_1,x_2)\in D_1 : \text{dist}(x,\partial D_1)=\text{dist}(x, E_1)\}.
\end{equation}
If we set $u^{(\eps)}(x):=u(x/\sqrt{\eps})$ and $A\slash\sqrt{\eps}:=\{ x:\sqrt{\eps} x\in A\}$, by scaling it is immediately seen that
\begin{eqnarray}\label{scaling2D}
\dys \eps^{p-2}\int_{A}|Du^{(\eps)}|^ph^{2-p} dx \!\! & = & \!\!\int_{D^0_1\slash\sqrt{\eps}}|Du|^px_2^{2-p}+ \int_{(D^0_1)^c\slash\sqrt{\eps}}|Du|^p\left(\frac{1}{\sqrt{\varepsilon}}-\sqrt{x_1^2+x_2^2}\right)^{2-p} dx \nonumber \\
\\
& =: & \!\! I_1^{\eps}+I_2^{\eps}. \nonumber
\end{eqnarray}
Notice that at least formally $I^{\eps}_1$ tends to $\dys \int_{\R^2_+}|Du|^p y_2^{2-p}dy$ as $\eps\to 0$ and we will control $I_2^{\eps}$, under suitable assumptions. This is the object of Section \ref{sec_main}. \medskip

In view of this scaling property, we consider the {optimal profile} problem, introduced in the Section \ref{sec_profilep}; that is,
\begin{eqnarray}\label{profilep2D}
\gamma_{p}=\inf\left\{\int_{\R^2_+}|{D} u|^p x_2^{2-p} dx+\int_{\R}V(Tu)d\hs^1 : u \in L^1_{\text{loc}}(\R^2_+) :  \hskip 2cm \right. \nonumber \\ 
\left. \lim_{t\to-\infty}\! Tu(t)=\alpha',  \ \lim_{t\to+\infty}\! Tu(t)=\beta'\right \}
\end{eqnarray}
and determines the line tension on the limit energy $\Phi$.

\medskip
\subsection{Compactness of the traces}

We prove the pre-compactness of the traces of the equi-bounded sequences for $H_\eps$, using the trace embedding of  $W^{1,p}(D_1, h^{2-p})$ in $W^{2-3/p,p}(\partial D_1)$ and the following lemma, which is an adaptation of \cite[Lemma 4.1]{palatucci07}, using the estimations in \cite[Lemma 4.1]{garroni06}.

\begin{lemma}\label{lem_flatcomp}
Let $(\ue)$ be a sequence in $W^{1,p}(D_1, h^{2-p}; [-m,m])$ and let $J\subset E_1$ be an open interval. For every $\delta$ such that $0<\delta<(\beta'-\alpha')/2$, define
$$
A_\eps:=\{x\in E_1:\tue(x)\leq\alpha+\delta\} \ \text{and} \ B_\eps:=\{x\in E_1:\tue(x)\geq\beta'-\delta\}
$$
and set
\begin{equation}\label{pesi}
a_\eps:=\frac{|A_\eps\cap J|}{|J|} \ \text{and} \ b_\eps:=\frac{|B_\eps\cap J|}{|J|}.
\end{equation}
Then
\begin{eqnarray}\label{nonoptimallb2}
H_\eps(\ue, D_1, J) \!\!\!  & \geq &  \!\!\! \left(\frac{S_p(\beta-\alpha-2\delta)^{p}}{(2p-3)(p-2)|J|^{2(p-2)}}\left(\!1-\frac{1}{(1-a_{\eps})^{2(p-2)}}-\frac{1}{(1-b_{\eps})^{2(p-2)}}\!\right) 
-C_1\right)\!\eps^{p-2} \nonumber \\
\\
& & \!+C_{\delta},  \nonumber
\end{eqnarray}
where $S_p$, $C_1$ and $C_{\delta}$ are positive constants not depending on $\eps$.
\end{lemma}
\medskip
\noindent
\\ {\bf Proof.} By the weighted Sobolev embedding of $W^{1,p}(D_1, h^{2-p})$ in $W^{2-3/p,p}(\partial D_1)$ (see \cite[Theorem 2.11]{nekvinda}), we have that there exists a constant $S_p$ such that for every $u \in W^{1,p}(D_1, h^{2-p})$
$$
\dys
\|Tu\|_{W^{2-3/p,p}(\partial D_1)} \leq S_p\|u\|_{W^{1,p}(D_1, h^{2-p})}.
$$
It follows that there exists a constant (still denoted by $S_p$) such that
\begin{equation*}
\dys \int_{D_1}|{D} u|^p h^{2-p}dx\geq S_p\int\!\!\int_{J\times J}\frac{|Tu(t)-Tu(t')|^p}{|t-t'|^{2(p-1)}}dt'dt-\int_{D_1}|u|^p h^{2-p} dx, \ \ \ \forall u \in W^{1,p}(D_1, h^{2-p}).
\end{equation*}
 Hence
\begin{eqnarray}\label{eq_nolb1}
\dys H_\eps(\ue, D_1, J)\! & = & \! \eps^{p-2}\int_{D_1}|{D}\ue|^p h^{2-p} dx+\frac{1}{\sqrt{\eps}}\int_{J}V(\tue)d\hs^1 \nonumber \\
\nonumber \\
&\geq& \! S_p\eps^{p-2}\int\!\!\int_{J\times J}\frac{|\tue(t)-\tue(t')|^p}{|t-t'|^{2(p-1)}}dt'dt+\frac{1}{\sqrt{\eps}}\int_{J}V(\tue)d\hs^1 \nonumber\\
\\ 
&  & \! -\eps^{p-2}\int_{D_1}|\ue|^ph^{2-p}dx \nonumber \\
\nonumber \\
& \geq & \! S_p\eps^{p-2}\int\!\!\int_{J\times J}\frac{|\tue(t)-\tue(t')|^p}{|t-t'|^{2(p-1)}}dt'dt+\frac{1}{\sqrt{\eps}}\int_{J}V(\tue)d\hs^1 -C_1\eps^{p-2}. \nonumber
\end{eqnarray}

The remaining part of the proof follows as in \cite[Lemma 4.1]{palatucci07}. \hfill $\Box$
\medskip

We are now in position to prove the compactness result stated in the following proposition.
\begin{proposition}\label{pro_comp_disco}
 If $(\ue)\subset W^{1,p}(D_1, h^{2-p};[-m,m])$ is a sequence such that $H_\eps(\ue)$ is bounded then $(\tue)$ is pre-compact in $L^1(E_1)$ and every cluster point belongs to $BV(E_1,\{\alpha',\beta'\})$.
\end{proposition}
\medskip
{\noindent
\bf Proof.}
By hypothesis, there exists a constant $C$ such that 
$ H_\eps(\ue)\leq C.$
 In particular
\begin{equation*}
\dys \int_{E_1}\!V(\tue)d\hs^1 \leq C\sqrt{\eps}
\end{equation*}
and this implies that
\begin{equation}\label{wconverge2}
V(\tue)\rightarrow 0 \ \text{in} \ L^{1}(E_1).
\end{equation}

Thanks to the growth assumptions on $V$, $(\tue)$ is equi-integrable. Hence, by Dunford-Pettis' Theorem, $(\tue)$ is weakly relatively compact in $L^{1}(E_1)$; i.e., there exists $v\in L^{1}(E_1)$ such that (up to subsequences) $\tue\rightharpoonup v$ in $ L^{1}(E_1)$.

 We have to prove that this convergence is strong in $L^{1}(E_1)$ and that $v \in BV(E_1; \{\alpha', \beta'\})$. This proof is standard, involving Young measures associated to sequences (see also \cite[Th\'eor\`eme 1-(i)]{alberti}). 
 Let $\nu_{x}$ be the Young measure associated with $(\tue)$.
 Since $V$ is a non negative continuous function in $\R$, 
we have
\begin{equation*}
\int_{E_1}\!\int_{\R}\!V(t)d\nu_{x}(t) \leq \liminf_{\eps \to 0}\int_{E_1}\!V(\tue)dx
\end{equation*} (see \cite[Theorem I.16]{valadier}).

Hence, by (\ref{wconverge2}), it follows that
$$
\int_{\R}V(t)d\nu_{x}(t)=0, \ \ \ \text{a.e.} \ x\in E_1,
$$
which implies the existence of a function $\theta$ on $[0,1]$ such that 
$$
\dys \nu_{x}(dt)=\theta(x)\delta_{\alpha'}(dt)+(1-\theta(x))\delta_{\beta'}(dt), \ \ \ x\in E_1 $${and}
$$
v(x)=\theta(x)\alpha'+(1-\theta(x))\beta', \ \ \ x \in E_1.
$$

It remains to prove that $\theta$ belongs to $BV(E_1;\{0,1\}).$
Let us consider the set $S$ of the points where approximate limits of $\theta$ is neither 0 nor 1.
For every $N\leq\hs^{0}(S)$ we can find $N$ disjoint intervals $\{J_{n}\}_{n=1,...,N}$ such that $J_{n}\cap S\neq\emptyset$ and such that the quantities $a_{\eps}^n$ and $b_{\eps}^n$, defined by (\ref{pesi}) replacing $J$ by $J_n$, satisfy
$$
a_{\eps}^n\to a^n \in (0,1) \ \ \ \text{and} \ \ \ b_{\eps}^n \to b^n \in (0,1) \ \ \
 \text{as}\ \eps \ \text{goes to zero.}
$$
We can now apply Lemma \ref{lem_flatcomp} in the interval $J_n$ and, taking the limit as $\eps\to0$ in the inequality (\ref{nonoptimallb2}), we obtain
$$
\dys \liminf_{\eps\to0} H_{\eps}(u_{\eps}, D_1, J_{n})\geq C_{\delta}.
$$
Finally, we use the sub-additivity of 
the non local part of the functional and we get 
\begin{equation*}
\liminf_{\eps\to0}H_{\eps}(\ue,D_1,E_1)\geq \sum_{n=1}^{N}\liminf H_{\eps}(\ue, D_1, J_{n}) \geq N C_{\delta}.
\end{equation*}
Since $(\ue)$ has equi-bounded energy, this implies that $S$ is a finite set. Hence, $\theta \in BV(E_1;\{0,1\})$ and the proof of the compactness for $H_{\eps}$ is complete. \hfill $\Box$

\medskip
\subsection{Lower bound inequality}

We will prove an optimal lower bound for $H_\eps$.
\begin{proposition}\label{pro_lb_disco}
For every $(u,v)$ in $BV(D_1; \{\al,\be\})\times BV(E_1; \{\al',\be'\})$ and every sequence $(\ue)\subset W^{1,p}(D_1, h^{2-p}; [-m,m])$ such that $\ue\to u$ in $L^1(D_1)$ and $\tue\to v$ in $L^1(E_1)$
\begin{equation}\label{lbdisco}
\dys \liminf_{\eps\to 0}H_\eps(\ue)\geq \gamma_p \hs^0(Sv).
\end{equation}
\end{proposition}
\noindent
\medskip
{\bf Proof.}
We will prove the lower bound inequality (\ref{lbdisco}) for $v$ such that
$$
v(t)=\begin{cases} \alpha', & \text{if} \ t\in (-1,0], \\ \beta', & \text{if} \ t \in (0,1).\end{cases} 
$$
Consider the natural extension of $v$ to the whole real line $\R$, still denoted by $v$; that is
$$
v(t)=\begin{cases} \alpha', & \text{if} \ t\leq 0, \\ \beta', & \text{if} \ t>0.\end{cases}
$$
\medskip

\noindent
{\it Step 0: Strategy of the proof.} We are looking for an extension of $\ue$ to the whole half-plane $\R^2_+$, namely $\we$, such that $\we$ is a competitor for $\gamma_p$ and $H_\eps(\we,\R^2_+,\R)\simeq  H_\eps(\ue, D_1, E_1)$ as $\eps\to0$ in a precise sense.
For every $\eps>0$, we will be able to find $s<1$ such that, for any given $\delta >0$ there exists $\eps_\delta>0$ and we have
\begin{eqnarray*}
H_\eps(\ue)\! &\geq& H_\eps(\ue, D_s, E_s) \\
&\geq&\! \gamma_p-\delta, \ \ \ \forall \eps\leq \eps_\delta.
\end{eqnarray*}
\medskip
\noindent
\\ {\em Step 1: Construction of the competitor.}
For every $s>0$, we denote by $\bar{u}$ the following extension of $v$ from $\R \setminus E_s$ to $\R^2_+\setminus D_s$ in polar coordinates
$$
\bar{u}(\rho,\theta):=\frac{\theta}{\pi^2}\alpha'+\left(1-\frac{\theta}{\pi^2}\right)\beta', \ \ \forall \theta \in [0,\pi), \ \forall \rho\geq s.
$$
We construct the competitor $\we$ simply gluing the function $\bar{u}$ and the function $\ue$. Hence, consider the cut-off function $\varphi$ in $C^{\infty}(\R^2_{+})$, such that $\varphi\equiv 1$ in ${D_s}$,  $\varphi\equiv 0$ in ${\R^2_+\setminus D_{s(\eps)}}$ and $|{D}\varphi|\leq 1/\eps^{\frac{p-2}{2(p-1)}}$, where we denote 
$$
\dys s(\eps):=s+\eps^{\frac{(p-2)}{2(p-1)}}.
$$
\begin{figure}[htp]
\centering
\psfrag{s1}{$\!-s\!\!-\!\!\eps^{\frac{p-2}{2(p-1)}}$} \psfrag{s2}{$-s$} \psfrag{s3}{$s$}  \psfrag{s4}{$\!\!\!s\!\!+\!\!\eps^{\frac{p-2}{2(p-1)}}$}
\psfrag{ue}{${u_\epsilon}$} \psfrag{ub}{$\bar{u}$} \psfrag{cc}{$\varphi\ue+(1-\varphi)\bar{u}$}
\includegraphics[scale=0.36]{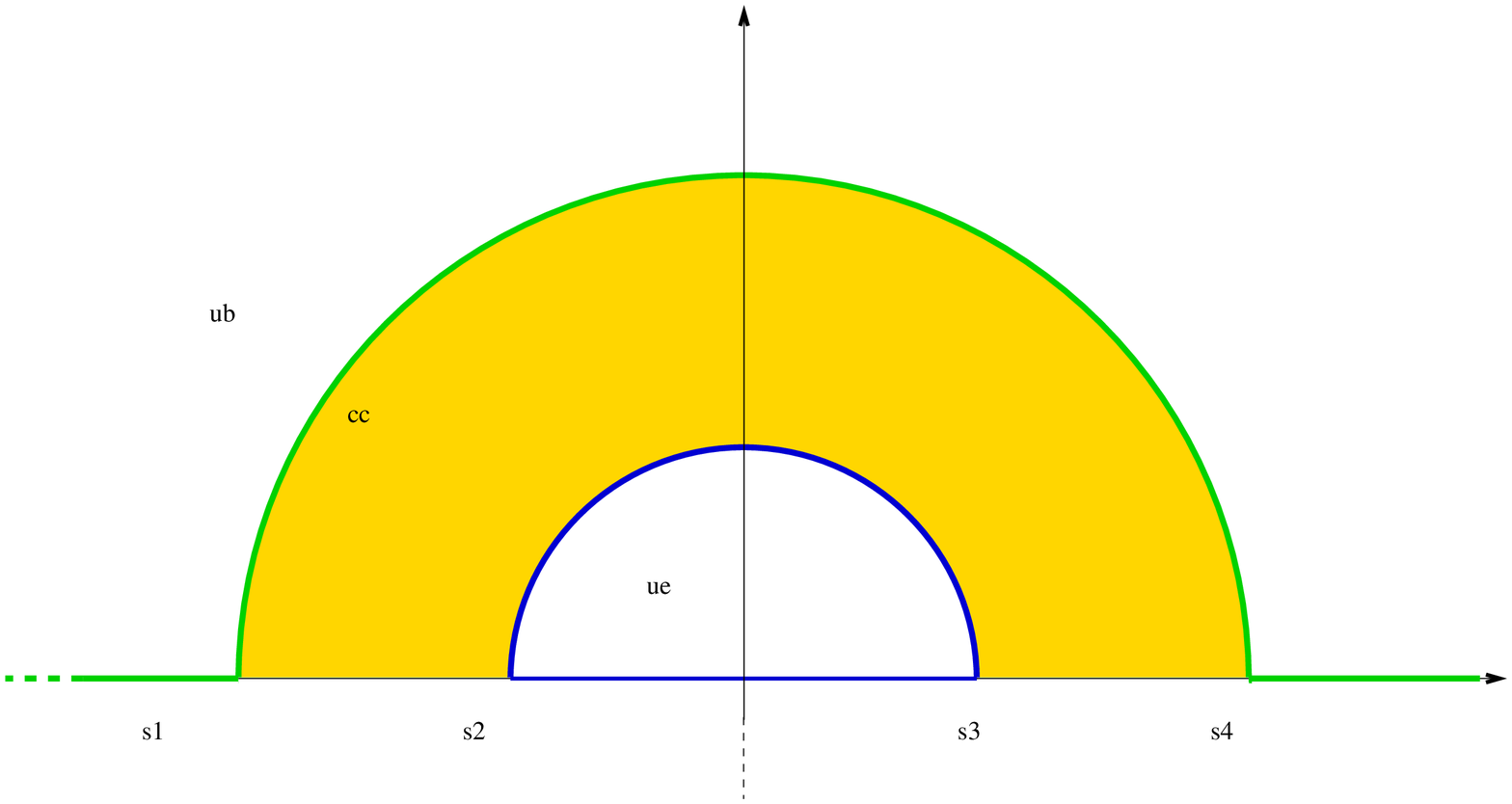}
\caption{The competitor $w_\eps$ (see also \cite[Fig. 2]{palatucci07}).}\label{fig_competitor2D}
\end{figure}

Thus, we consider
$$
\dys \we:=\begin{cases}
\ue &\text{in} \ D_s, \\
\varphi\ue+(1-\varphi)\bar{u} & \text{in} \ D_{s(\eps)}\setminus \ D_{s},\\
\bar{u} &\text{in} \ \R^2_+\setminus D_{s(\eps)}.
\end{cases}
$$
Note that  $\dys \lim_{t\to-\infty}\twe(t)=\alpha'$ and  $\dys \lim_{t\to+\infty}\twe(t)=\beta'$. \medskip
\noindent
\\ {\em Step 2: Choice of the annulus.} We need to choose an annulus in the half-disk, in which we can recover a suitable quantity of energy of $\ue$. Since $\ue$ has equi-bounded energy $H_\eps(\ue)$ in $D_1$, there exists $L>0$ such that $\forall \eps>0$  $\exists s\in \dys \left(\frac{1}{2}, 1-\eps^{\frac{p-2}{2(p-1)}}\right)$ such that
\begin{equation}\label{eq_step2}
\dys \eps^{p-2}\int_{D_{s(\eps)}\setminus D_s}|{D} \ue|^p h^{2-p} dx+\frac{1}{\sqrt{\eps}}\int_{E_{s(\eps)}\setminus E_s}V(T\ue)d\hs^1\leq L\eps^{\frac{p-2}{2(p-1)}}.
\end{equation}
\medskip
\noindent
\\ {\em Step 3: Estimates.}
For every $s$ like in Step 2, we have
\begin{equation}\label{eq_step3-dist}
H_\eps(\ue)  \geq  \eps^{p-2}\int_{D_s\cap D_1^{0}}|D\ue|^p h^{2-p} dx+\frac{1}{\sqrt{\eps}}\int_{E_s}V(\tue)d\hs^1
\end{equation}
where $D^0_1$ is defined by (\ref{def_d10}).
\smallskip

By the scaling property of $H_\eps$ (see (\ref{scaling2D})) and noticing that here we have always $h^{2-p}=x_2^{2-p}$, we have
\begin{eqnarray}\label{eq_step3-pro}
\eps^{p-2}\!\int_{D_s\cap D_1^{0}}\!|D\ue|^p x_2^{2-p} dx+& \!\!\! \dys \frac{1}{\sqrt{\eps}}& \!\!\!\int_{E_s}\!V(\tue)d\hs^1 \nonumber \\
& = & \!\! \int_{(D_s\cap D_1^{0})/\sqrt{\eps}}|D\ue^{(\eps)}|^p y_2^{2-p} dy+\int_{E_s}V(\twe^{(\eps)})d\hs^1 \nonumber\\
\nonumber \\
& = & \!\! H_1(w_\eps^{(\eps)}, \R^2_+,\R) - \int_{(\R^2_+\setminus(D_s\cap D_1^0))/\sqrt{\eps}}|D\we^{(\eps)}|^p y_2^{2-p} dy \\
& & \!\! -\int_{(\R\setminus E_s)/{\sqrt{\eps}}}V(\we^{(\eps)})d\hs^1 \nonumber\\
\nonumber\\
& \geq & \!\! \gamma_p-\eps^{p-2}\!\int_{\R^2_+\setminus(D_s\cap D_1^0)}\!\!|D\we|^p x_2^{2-p} dx 
-\frac{1}{\sqrt{\eps}}\!\int_{\R\setminus E_s}\!\!V(\twe)d\hs^1,\nonumber
\end{eqnarray}
where we recall that $u^{(\eps)}$ is the rescaled function defines by $u^{(\eps)}=u(x/\sqrt{\eps})$.
\smallskip

Using (\ref{eq_step3-dist}) and (\ref{eq_step3-pro}), we get
\begin{eqnarray*}
\dys H_\eps(\ue) \!\! & \geq & \!\! \gamma_p-\eps^{p-2}\int_{\R^2_+\setminus(D_s\cap D_1^0)}|D\we|^p x_2^{2-p} dx -\frac{1}{\sqrt{\eps}}\int_{\R\setminus E_s}V(\twe)d\hs^1 \nonumber\\ \nonumber\\
& \geq & \!\! \gamma_p -\eps^{p-2}\int_{\R^2_+\setminus (D_{s(\eps)}\cap D_1^0)}|D\bar{u}|^p x_2^{2-p} dx -\eps^{p-2}\int_{(D_{s(\eps)}\setminus D_{s})\cap D_1^0}|D\we|^p x_2^{2-p} dx \nonumber \\\\
& & \!\! -\frac{1}{\sqrt{\eps}}\int_{\R\setminus E_{s}}V(T\we) d\hs^1. \nonumber
\end{eqnarray*}
\smallskip

Thus
\begin{eqnarray}\label{eq_step3-gamma}
\dys \gamma_p \!\! & \leq & \!\! H_\eps(\ue) + \eps^{p-2}\int_{(\R^2_+\setminus D_{s(\eps)})\cap D_1^0}|D\bar{u}|^p x_2^{2-p}dx  \nonumber \\ \\
 & \!\! & + \eps^{p-2}\int_{(D_{s(\eps)}\setminus D_{s})\cap D_1^0}|D\we|^p x_2^{2-p} dx
+ \frac{1}{\sqrt{\eps}}\int_{\R\setminus E_{s}}V(T\we) d\hs^1. \nonumber
\end{eqnarray}
\smallskip

Using the fact that $(X+Y+Z)^p\leq 3^{p-1}(X^p+Y^p+Z^p)$ and the definition of $\we$, the second integral in the right hand side of (\ref{eq_step3-gamma}) can be estimated as follows
\begin{eqnarray*}
\dys \eps^{p-2}\int_{(D_{s(\eps)}\setminus D_{s})\cap D_1^0}|D\we|^p x_2^{2-p} dx  \!\! & \leq & \!\! 3^{p-1}\eps^{p-2}\int_{(D_{s(\eps)}\setminus D_{s})\cap D_1^0}|D\bar{u}|^p x_2^{2-p} dx \\
& & \!\! + 3^{p-1}\eps^{p-2}\int_{(D_{s(\eps)}\setminus D_{s})\cap D_1^0}|D\ue|^p x_2^{2-p} dx \\\\
& & \!\! + 3^{p-1}C_1(s+\eps^{\frac{p-2}{2(p-1)}})^{4-p}\eps^{\frac{(p-2)^2}{2(p-1)}}.
\end{eqnarray*}
Here we also used that $h(x_1,x_2)^{2-p}=x_2^{2-p}$ with $x_2=\rho\sin{\theta}$, $\frac{s}{1+\sin{\theta}}\leq \rho \leq \frac{s(\eps)}{1+\sin{\theta}}$ and it gives the estimate
\begin{eqnarray*}
\dys \int_{(D_{s(\eps)}\setminus D_{s})\cap D_1^0}|D\varphi|^p|\ue-\bar{u}|^p h^{2-p} dx \!\! 
& \leq & \!\! \frac{(2m)^p}{\eps^{\frac{p(p-2)}{2(p-1)}}}\int_0^\pi\!\!\int_{s}^{s(\eps)}\rho^{2-p}\sin^{2-p}{\theta}\rho d\rho d\theta \\
\nonumber \\
& \leq & \! \frac{(2m)^p s(\eps)^{4-p}}{(4-p)\eps^{\frac{p(p-2)}{2(p-1)}}}\int_0^{\pi}\frac{1}{\sin^{p-2}(\theta)}d\theta \\
\\
& =  & \!\! \frac{C_1(s+\eps^{\frac{p-2}{2(p-1)}})^{4-p}}{\eps^{\frac{p(p-2)}{2(p-1)}}}.
\end{eqnarray*}

\smallskip

It follows
\begin{eqnarray}\label{eq_step3-ubar}
\dys \gamma_p \!\! & \leq & \!\!  H_\eps(\ue)+3^{p-1}\eps^{p-2}\int_{\R^2_+\setminus D_s \cap D_1^0}|D\bar{u}|^p x_2^{2-p}dx
 +3^{p-1}\eps^{p-2}\int_{(D_{s(\eps)}\setminus D_{s})\cap D_1^0}|D\ue|^p x_2^{2-p} dx \nonumber \\\\
& &  \!\! +\frac{1}{\sqrt{\eps}}\int_{\R\setminus E_{s}}V(T\we) d\hs^1 + 3^{p-1}C_1(s+\eps^{\frac{p-2}{2(p-1)}})^{4-p}\eps^{\frac{(p-2)^2}{2(p-1)}}. \nonumber
\end{eqnarray} \smallskip
Using the definition of $\bar{u}$, we can compute the first integral
\begin{eqnarray*}
\dys \int_{\R^2_+\setminus D_s\cap D_1^0}|D\bar{u}|^p x_2^{2-p} dx \! \! & = & \!\! \frac{C_2}{s^{2(p-2)}}
\end{eqnarray*}
and (\ref{eq_step3-ubar}) becomes
\begin{eqnarray}\label{eq_step3-ubar2}
\dys \gamma_p \!\! & \leq & \!\!  H_\eps(\ue)+3^{p-1}\frac{C_2}{s^{2(p-2)}}
+3^{p-1}\eps^{p-2}\int_{(D_{s(\eps)}\setminus D_{s})\cap D_1^0}|D\ue|^p x_2^{2-p} dx \nonumber \\
\\
& & \!\!+\frac{1}{\sqrt{\eps}}\int_{E_{s(\eps)}\setminus E_{s}}V(T\we) d\hs^1
+ 3^{p-1}C_1\eps^{\frac{(p-2)^2}{2(p-1)}}. \nonumber
\end{eqnarray}

Let us estimate the second integral in the right hand side of (\ref{eq_step3-ubar2}).
Since $\twe=\al'$ and $\twe=\be'$ on $\R\setminus E_{s(\eps)}$, we have
$$
\dys \frac{1}{\sqrt{\eps}}\int_{\R\setminus E_{s(\eps)}}V(T\bar{u})d\hs^1+\frac{1}{\sqrt{\eps}}\int_{E_{s(\eps)}\setminus E_s}V(\twe)d\hs^1\equiv \frac{1}{\sqrt{\eps}}\int_{E_{s(\eps)}\setminus E_s}V(\twe)d\hs^1.
$$\smallskip
For every $\delta>0$, let us define
$$
\dys E_\delta:=\left\{x\in E_{s(\eps)}\setminus E_s:|\tue-\beta'|>\delta \ \text{and} \ |\tue-\alpha'|>\delta\right\}.
$$
Thanks to Step 2,
there exists $\dys N>\frac{L}{\omega_\delta}$ $\bigl(\text{where we denote by} \ \dys \omega_\delta:=\!\!\min_{\stackrel{|t-\al'|\geq\delta} {|t-\be'|\geq\delta}}V(t)\bigl)$  such that $\forall \delta>0$ $\exists \eps_\delta$ such that
\begin{equation}\label{eq_lemout}
\dys |E_\delta|\leq N\eps^{\frac{p-2}{2(p-1)}}\sqrt{\eps}, \ \ \ \forall\eps\leq\eps_\delta. 
\end{equation}

In particular, choosing $\delta$ small, the convexity of $V$ near its wells provides
\begin{eqnarray}\label{eq_step3i2}
\!\frac{1}{\sqrt{\eps}}\int_{(E_{s(\eps)}\setminus E_s)\setminus E_\delta}\!\!V(\twe)d\hs^1+\frac{1}{\sqrt{\eps}}\int_{E_{\delta}}\!V(\twe)d\hs^1  \! & \leq & \! \frac{1}{\sqrt{\eps}}\int_{(E_{s(\eps)}\setminus E_s)\setminus E_\delta}\!\!V(\tue)d\hs^1 \nonumber \\
\\
& & \! +\omega_m N\eps^{\frac{p-2}{2(p-1)}},\nonumber
\end{eqnarray}
where $\dys \omega_m\!:=\!\max_{|t|<m}V(t)$ and we used the inequality (\ref{eq_lemout}).

\medskip

Finally, by (\ref{eq_step2}), (\ref{eq_step3-ubar2}) and (\ref{eq_step3i2}), we obtain, for every $\delta >0$
\begin{eqnarray*}
\dys H_\eps(\ue) \! &\geq &\! \gamma_p -   \left(3^{p-1}\left(\eps^{p-2}\int_{D_{s(\eps)}\setminus D_s}\!|{D} \ue|^p h^{2-p} dx+\frac{1}{\sqrt{\eps}}\int_{E_{s(\eps)}\setminus E_s}\!V(T\ue)d\hs^1\right)
\hskip 2cm \right. \nonumber \\
\nonumber \\ 
& & \! \dys \left. {3^{p-1}}\frac{C_2}{s^{2(p-2)}}\eps^{p-2} + {3^{p-1}}C_1(s+\eps^{\frac{p-2}{2(p-1)}})^{4-p}\eps^{\frac{(p-2)^2}{2(p-1)}} + \omega_m N\eps^{\frac{p-2}{2(p-1)}}\right) \nonumber \\
& \geq & \! \dys \gamma_p -  \left(3^{p-1}L\eps^{\frac{p-2}{2(p-1)}}
+  {3^{p-1}}\frac{C_2}{s^{2(p-2)}}\eps^{p-2} + {3^{p-1}}C_1(s+\eps^{\frac{p-2}{2(p-1)}})^{4-p}\eps^{\frac{(p-2)^2}{2(p-1)}} \hskip 2cm \right. \nonumber \\
\\
& & \! \left. + \omega_m N\eps^{\frac{p-2}{2(p-1)}}\right). \nonumber
\end{eqnarray*}

Notice that for every $\eps>0$, $\dys s\in\left({1}/{2},1-\eps^{\frac{p-2}{2(p-1)}}\right)$. Hence, taking the limit as $\eps\to0$, we get
$\dys \liminf_{\eps \to 0} H_{\eps}(\ue) \geq\gamma{_p},$ which concludes the proof in the case of a function $v$ with one jump, i.e. $\mathcal H^0(S v)=1$. The case $\mathcal H^0(S v)>1$ can be treated similarly.  \hfill $\Box$

\medskip
\subsection{Reduction to the flat case}

According to the idea of Alberti, Bouchitte and Seppecher in \cite{alberti98}, we will prove the Theorem \ref{c4_teoremap} after arguments of slicing and blow-up. More precisely, it is possible to deform each neighborhood in a bi-Lipschitz fashion in order to straighten the boundary of $\Om$, without changing much the functional (see \cite[Proposition 4.9]{alberti98} and \cite[Proposition 5.2]{palatucci07}). 

To this aim, we recall the definition of the ``isometry defect'', introduced by Alberti, Bouchitt\'e and Seppecher\cite{alberti98}.
\smallskip

As usual, we denote by  $O(3)$ the set of linear isometries on $\R^3$.

\begin{definition}\label{def_isometry} Let $A_1, A_2 \subset \R^3$ and let $\Psi:\overline{A_1}\to\overline{A_2}$ bi-Lipschitz homeomorphism. Then the {\em ``isometry defect  $\delta({\Psi})$ of $\Psi$''} is the smallest constant $\delta$ such that
\begin{equation}\label{eq_422abs}
{\text{\em dist}}({D}\Psi(x), O(3))\leq\delta, \ \ \text{for} \ \text{a.e.} \ x\in A_1. 
\end{equation}
\end{definition}
Here ${D}\Psi(x)$ is regarded as a linear mapping of $\R^3$ into $\R^3$. The distance between linear mappings is induced by the norm $\|\cdot\|$, which, for every $L$, is defined as the supremum of $|Lv|$ over all $v$ such that $|v| \leq 1$. Hence, for every $L_1, L_2:\R^3\to\R^3$:
$$
\dys \text{dist}(L_1,L_2):=\dys \sup_{x:|x|\leq 1}|L_1(x)-L_2(x)|.
$$

By (\ref{eq_422abs}), we get
\begin{equation}\label{eq_iso1}
\|{D}\Psi(x)\|\leq 1+\delta(\Psi) \  \text{for a.e.}\ x \in A_1,
\end{equation}
and then $\Psi$ is $(1+\delta(\Psi))$-Lipschitz continuous on every convex subset of $A_1$. Similarly, $\Psi^{-1}$ is $(1-\delta(\Psi))^{-1}$-Lipschitz continuous on every convex subset of $A_2$.
\smallskip


The following proposition shows that the localized energy $F_\eps(u,B_r(x)\cap\Om, B_r(x)\cap \partial\Om)$ can be replaced by the energy $F_\eps(u, D_r, E_r)$, where $B_r$ is the two-dimensional ball of radius $r$ centered in the origin.

\begin{proposition}\label{pro_abs410}
For every $x \in \partial\Om$ and every positive $r$ smaller than a certain critical value $r_x>0$, there exists a bi-Lipschitz map $\Psi_r:\overline{{D}_r}\to\overline{{\Om\cap B_r(x)}}$ such that
\begin{itemize}
\item[(a)]{$\Psi_r$ takes ${D}_r$ onto $\Om\cap B_r(x)$ and ${E}_r$ onto $\partial\Om\cap B_r(x)$;}
\item[(b)]{$\Psi_r$ is of class $C^1$ on $D_r$ and $\|D\Psi_r-I\|\leq\delta_r$ everywhere in $D_r$, where $\delta_r\to0$ as $r\to0$.}
\end{itemize}
In particular, the isometry defect of $\Psi_r$ vanishes as $r\to 0$. Moreover,
$$
F_\eps(u,B_r(x)\cap\Om, B_r(x)\cap \partial\Om)\geq (1-\delta(\Psi))^{p+3}F_\eps(u\circ\Psi, D_r, E_r).
$$
\end{proposition}
The proof is a simple modification of the one by Alberti, Bouchitt\'e and Seppecher in \cite{alberti98}, Proposition 4.9 and Proposition 4.10, where they treat the case $p=2$ (see also \cite[Proposition 6.1]{gonzalez}).

\medskip

Finally, we need to prove compactness and a lower bound inequality for the following energies 

$$
\dys F_\eps(u, \mathcal{D}, \mathcal{E})=\eps^{p-2}\int_{\mathcal{D}}|Du|^p h^{2-p} dx+\frac{1}{\eps^{\frac{p-2}{p-1}}}\int_{\mathcal{D}}W(u) h^ {\frac{p-2}{p-1}} dx+\frac{1}{\sqrt{\eps}}\int_{\mathcal{E}}V(Tu)d\hs^2,
$$
where $\mathcal{D}\subset\R^3$ is the open half-ball centered in 0 with radius $r>0$ and $\mathcal{E}\subset\R^2$ is defined by
$$
\dys \mathcal{E}:=\left\{(x_1,x_2,x_3)\in\R^3 : |x|\leq r, x_3=0\right\}.
$$
We will reduce to Proposition \ref{pro_comp_disco} and Proposition \ref{pro_lb_disco} via a suitable slicing argument. 
\begin{figure}[htp!]
\centering
\psfrag{dset}{$\mathcal{D}$} \psfrag{eset}{$\mathcal{E}$} \psfrag{dy}{$\mathcal{D}^y$} \psfrag{y}{$y$} \psfrag{x3}{\!$x_3$} \psfrag{ee}{$\mathcal{E}_e$} \psfrag{ey}{$\mathcal{E}^y$} \psfrag{m}{\!\!\!$M$}\psfrag{e}{\!$e$}
\includegraphics[scale=0.16]{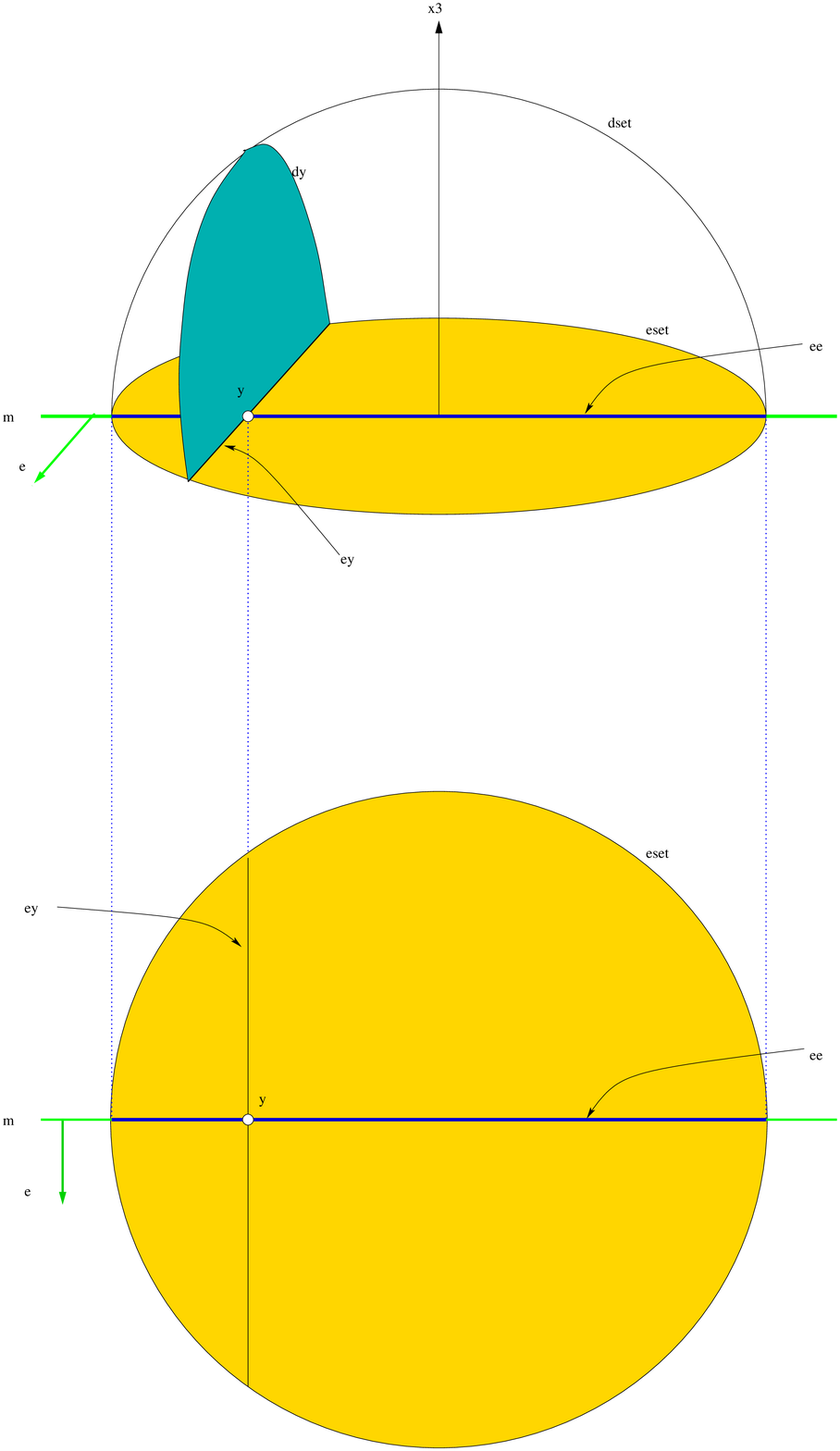}
\caption{The sets $\mathcal{D}, \mathcal{E}, \mathcal{E}_e, \mathcal{E}^y$ and $\mathcal{D}^y$ (see also \cite[Fig. 4]{alberti98}).}\label{fig_slice3d}
\end{figure}
\noindent
\\ We use the following notation: $e$ is a unit vector in the plane $P:=\{x_3=0\}$; $M$ is the orthogonal complement of $e$ in $P$; $\pi$ is the projection of $\R^3$ onto $M$; for every $y\in \mathcal{E}_e:=\pi(\mathcal{E})$, we denote by $\mathcal{E}^y:=\pi^{-1}(y)\cap \mathcal{E}, \ \mathcal{D}^y:=\pi^{-1}(y)\cap \mathcal{D}$ (see Fig. \ref{fig_slice3d}); for every function $u$ defined on $\mathcal{D}$ we consider the trace of $u$ on $\mathcal{E}^y$, i.e., the one-dimensional function
$$
u^y_e(t):=u(y+te).
$$

\begin{proposition}\label{pro_slice}
Let $(\ue)\subset W^{1,p}(\mathcal{D},h^{2-p};[-m,m])$ be a sequence with uniformly bounded energies $F_\eps(\ue, \mathcal{D}, \mathcal{E})$. Then the traces $\tue$ are pre-compact in $L^1(\mathcal{E})$ and every cluster point belongs to $BV(\mathcal{E}; \{\al',\be'\})$. Moreover, if $\tue\to v$ in $L^1(\mathcal{E})$, then
\begin{equation}\label{ineq_slice}
\dys \liminf_{\eps\to 0}F_\eps(\ue, \mathcal{D}, \mathcal{E})\geq\gamma_p\left|\int_{\mathcal{E}\cap Sv}\nu_v\right|d\hs^1.
\end{equation}
\end{proposition}
\medskip
\noindent {\bf Proof.}
By Fubini's Theorem, for every $\eps>0$, we get
\begin{eqnarray}\label{slice_step1}
\dys F_\eps(\ue, \mathcal{D}, \mathcal{E}) & \geq & \eps^{p-2}\int_\mathcal{D}|{D} \ue|^p h^{2-p} dx+\frac{1}{\sqrt{\eps}}\int_{\mathcal{E}} V(\tue)d\hs^2 \nonumber \\
\nonumber \\
\dys & \geq & \int_{\mathcal{E}_e}\left[\eps^{p-2}\int_{\mathcal{D}^y}|{D} u_\eps^y|^p h^{2-p} dx+\frac{1}{\sqrt{\eps}}\int_{\mathcal{E}^y} V(Tu_\eps^y)d\hs^1\right]dy \nonumber\\
\\
\dys & = & \int_{\mathcal{E}_e}H_\eps(\ue^y, \mathcal{D}^y, \mathcal{E}^y)dy\nonumber
\end{eqnarray}
and the 2D functional in the integration above has been studied in the last section.
The remaining part of the proof follows exactly as in \cite[Proposition 4.7]{alberti98}.\hfill $\Box$

\bigskip

\subsection{Existence of an optimal profile}

We conclude this section with the proof of the existence of a minimum for the optimal profile problem  (\ref{profilep2D}), showing that the minimum for $\gamma_p$ is achieved by a function with non-decreasing trace.

\begin{proposition}\label{pro_min2d}
The minimum for $\gamma_p$ defined by (\ref{profilep2D}) is achieved by a function $u$ such that $Tu$ is a non-decreasing function in $\R$.
\end{proposition}
\medskip
\noindent
\\ {\bf Proof.}  Note that, since the energy $H_1$ is decreasing under truncation by $\al'$ and $\be'$, it is not restrictive to minimize the problem (\ref{profilep2D}) with the additional condition $\dys \al'\leq u\leq\be'$.

We denote by 
$$
\dys X:=\left\{w:\R^2_+\to[\al',\be'] : w\in L^1_{\text{loc}}(\R^2_+),  \lim_{t\to-\infty}Tw(t)=\al',  \lim_{t\to+\infty}Tw(t)=\be' \right\}
$$
$$
\dys X^{\ast}:=\left\{w\in X : Tw \ \text{is non-decreasing}, \ Tw(t)\geq\frac{\al'+\be'}{2} \ \text{for} \ t>0, \ Tw(t)\leq\frac{\al'+\be'}{2} \ \text{for} \ t<0\right\}.
$$

\medskip

Let $u$ be in $X$, we denote by $u^{\star}$ its monotone increasing rearrangement in direction $x_1$. Since monotone increasing rearrangement in one direction decreases the weighted $L^p$-norm of the gradient (see \cite[Theorem 3]{berestycki}), the infimum of $H_1$ on $X$ is equal to the infimum of $H_1$ on $X^{\ast}$.

\medskip

Now we can prove by Direct Method that the infimum of $H_1$ on $X^{\ast}$ is achieved.
 
Take a minimizing sequence $(u_n)\subset X^{\ast}$.
In particular, $H_1(u_n,\R^2_+,\R)\leq C$, $Du_n$ converges weakly to $Du$ in $L^p(\R^2_+, h^{2-p})$ and $u_n$ converges to $u$ weakly in $W^{1,p}_{\text{loc}}(\R^2_+, h^{2-p})$.
Since $\dys \int_{\R^2_+}|Du_n|^p h^{2-p} dx$ is bounded, we can find a function $u\in L^{1}_{\text{loc}}(\R^2_+, h^{2-p})$
such that (up to a subsequence)
$$
\dys Du_n\weak Du \ \text{in} \ L^p(\R^2_+, h^{2-p}) \ \text{and} \ u_n \weak u \ \text{in} \ L^p_{\text{loc}}(\R^2_+, h^{2-p}).
$$
By the 
trace embedding of $W^{1,p}(\R^2_+, h^{2-p})$ in $W^{2-3/p,p}(\R)$, 
we have
$$
\dys Tu_n\weak Tu \ \text{in} \ W^{2-3/p,p}_{\text{loc}}(\R).
$$
By the compact embedding of  $C^{0}_{\text{loc}}(\R)$ in $W^{2-3/p,p}_{\text{loc}}(\R)$ (see \cite[Theorem 7.34]{adams}), 
we have that, up to a subsequence, $Tu_n$ locally uniformly converges to $Tu$. Thus $Tu$ is non-decreasing and satisfies
$$
\dys Tu(t)\geq\frac{\al'+\be'}{2} \ \text{for}  \ t>0 \ \ \ \text{and} \ \ \ \dys Tu(t)\leq\frac{\al'+\be'}{2} \ \text{for} \ t<0.
$$
Let us show that $\dys \lim_{t\to-\infty}Tu(t)=\al'$ and $\dys \lim_{t\to+\infty}Tu(t)=\be'$. Since $Tu$ is non-decreasing in $[\al',\be']$, there exist $\dys a\leq \frac{\al'+\be'}{2}$ and $\dys b\geq \frac{\al'+\be'}{2}$ such that
$$
\dys a:=\lim_{t\to-\infty}Tu(t) \ \ \text{and} \  \ b:=\lim_{t\to+\infty}Tu(t).
$$
By contradiction, we assume that either $a\neq\al'$ or $b\neq\be'$. Then, since $V$ is continuous and strictly positive in $(\al',\be')$, we obtain
$$
\dys \int_{\R}V(Tu)d\hs^1=+\infty,
$$
This is impossible, because, by Fatou's Lemma, we have
$$
\dys \int_{\R}V(Tu)d\hs^1\leq\liminf_{n\to+\infty} \int_{\R}V(Tu_n)d\hs^1< \liminf_{n\to+\infty}H_1(u_n,\R^2_+,\R) <+\infty.
$$
Hence, $u$ is in $X^{\ast}$. Since $H_1$ is clearly lower semicontinuous on sequences such that $Du_n\weak Du$ in $L^p(\R^2_+,h^{2-p})$ and $Tu_n\to Tu$ pointwise,
this concludes the proof.\hfill $\Box$

\bigskip

\section{Proof of the main result}\label{sec_main}

In the previous sections, we have obtained the main ingredients of the proof of Theorem \ref{c4_teoremap}, which follows as in the quadratic case in \cite{alberti98}, with the needed modifications due to the presence of the weight (like in \cite{gonzalez}) and to the super-quadratic growth in the singular perturbation term (like in \cite{palatucci07}).

\bigskip

\subsection{Compactness}

Let a sequence $(\ue)\subset W^{1,p}(\Om)$ be given such that $F_\eps(\ue)$ is bounded.
Since $F_\eps(\ue)\geq F_\eps(\ue,\Om,\emptyset)\equiv G_\eps(\ue,\Om)$, by the statement (i) of Theorem \ref{teo_modicap}, the sequence $(\ue)$ is pre-compact in $L^1(\Om)$ and there exists $u\in BV(\Om; \{\al,\be\})$ such that $\ue\to u$ in $L^1(\Om)$.

On the boundary, by slicing, we may use Proposition \ref{pro_slice}, that implies that $(\tue)$ is pre-compact in $L^1(\partial\Om)$ and that its cluster points are in $BV(\partial\Om; \{\al',\be'\})$.\hfill $\Box$

\bigskip

\subsection{Lower bound inequality}

The proof of the lower bound inequality of Theorem \ref{c4_teoremap} follows by putting together the results in the interior, Theorem \ref{teo_modicap}-(ii) and Proposition \ref{pro_modica2p}-(i), and the ones on the boundary, Section \ref{ch_2d} and Proposition \ref{pro_lb_disco}.
\medskip

Let a sequence $(\ue)\subset W^{1,p}(\Om, h^{2-p})$ be given such that $\ue\to u \in BV(\Om, \{\al,\be\})$ in $L^1(\Om)$ and $\tue\to v\in BV(\partial\Om, \{\al'\be'\})$ in $L^1(\partial\Om)$. We have to prove that
\begin{equation}\label{eq_lb}
\dys \liminf_{\eps\to 0}F_\eps(\ue)\geq\Phi(u,v),
\end{equation}
where $\Phi$ is given by (\ref{fu_phi}).
\\ Clearly, we can assume that $\dys\liminf_{\eps\to0}F_\eps(\ue)<+\infty$.
\medskip

For every $\eps>0$, let $\mu_\eps$ be the energy distribution associated with $F_\eps$ with configuration $\ue$; i.e., $\mu_\eps$ is the positive measure given by
\begin{equation}\label{distribution}
\dys \mu_\eps(B):=\eps^{p-2}\int_{\Om\cap B}|{D}\ue|^p dx+\frac{1}{\eps^{\frac{p-2}{p-1}}}\int_{\Om\cap B}W(\ue) dx+\frac{1}{\sqrt \eps}\int_{\partial\Om\cap B}V(\tue) d\hs^2,
\end{equation}
for every $B\subset\R^3$ Borel set.

Similarly, let us define
\begin{eqnarray*}
 \dys\mu^1(B)\!\!&:=&\!\!\sigma_p\hs^2(Su\cap B),\\
 \dys \mu^2(B)\!\!&:=&\!\! c_p\int_{\partial\Om\cap B}|\mathcal{W}(Tu)-\mathcal{W}(v)|d\hs^2,\\
 \dys \mu^3(B)\!\!&:=&\!\!\gamma_p\hs^1(Sv\cap B).
\end{eqnarray*}
The total variation $\|\mu_\eps\|$ of the measure $\mu_\eps$ is equal to $F_\eps(\ue)$, and $\|\mu^1\|+\|\mu^2\|+\|\mu^3\|$ is equal to $\Phi(u,v)$. The quantity 
$\|\mu_\eps\|$ is bounded and we can assume that $\mu_\eps$ converges in the sense of measures to some finite measure $\mu$. 
Then, by the lower semicontinuity of the total variation, we have
$$
\liminf_{\eps\to0}F_\eps(\ue)\equiv\liminf_{\eps\to0}\|\mu_\eps\|\geq\|\mu\|.
$$
Since the measures $\mu^i$ are mutually singular, we obtain the lower bound inequality (\ref{eq_lb}) if we prove that
\begin{equation}\label{ineq_distr}
\dys \mu\geq\mu^i, \ \ \text{for} \ i=1, 2, 3.
\end{equation}

It is enough to prove that $\mu(B)\geq\mu^i(B)$ for all sets $B\subset\R^3$ such that $B\cap\Om$ is a Lipschitz domain and $\mu(\partial B)=0$. 
\medskip

By the inequality of statement (ii) of Theorem \ref{teo_modicap}, we have
$$
\mu(B) = \lim_{\eps\to0}\mu_\eps(B)\geq\liminf_{\eps\to0} F_\eps(\ue, \Om\cap B, \emptyset)\geq\sigma_p\hs^2(Su\cap B)\equiv\mu^1(B).
$$\medskip

Similarly, we can prove that $\mu\geq\mu^2$. We have
\begin{eqnarray*}
\mu(B) & = & \lim_{\eps\to0}\mu_\eps(B) \ \geq \ \liminf_{\eps\to0} F_\eps(\ue, \Om\cap B, \emptyset) \\
& \geq & c_p\int_{\partial\Om\cap B}|\mathcal{W}(Tu)-\mathcal{W}(v)| d\hs^2 \ \equiv \ \mu^2(B),
\end{eqnarray*}
where we used Proposition \ref{pro_modica2p}-(i) with $A:=B\cap\Om$ and $A':=B\cap\partial\Om$.\medskip

The inequality $\mu\geq\mu^3$ requires a different argument. Notice that $\mu^3$ is the restriction of $\hs^1$ to the set $Sv$, multiplied by the factor $\gamma_p$.
Thus, if we prove that
\begin{equation}\label{ineq_mu3}
\dys\liminf_{r\to0}\frac{\mu(B_r(x))}{2r}\geq\gamma_p, \ \ \hs^1\text{-a.e.} \ x\in Sv,
\end{equation}
for $B_r(x)$ as in Proposition \ref{pro_abs410}, we obtain the required inequality.

Let us fix $x\in Sv$ such that there exists $\dys\lim_{r\to0}\frac{\mu(B_r(x))}{2r}$ and $Sv$ has one-dimensional density equal to $1$. We denote by
$\nu_v$ the unit normal at $x$.

For $r$ small enough, we choose a map $\Psi_r$ such as in Proposition \ref{pro_abs410}. Thus we have $\dys \Psi_r(\overline{D_r})={\Om\cap B_r(x)}, \ \Psi_r(E_r)=\partial\Om\cap B_r(x)$ and $\delta(\Psi_r)\to0$ as $r\to0$.

Let us set
$$
\dys {\bar{u}_{\eps}}:=\ue \circ \Psi_r \  \text{and}  \ \bar{v}:=v \circ \Psi_r.
$$
Hence, $T{\bar{u}_{\eps}}\to\bar{v}$ in $L^1(E_r)$ and $\bar{v}\in BV(E_r, \{\al',\be'\})$.
So, thanks to Proposition \ref{pro_abs410}, we obtain
\begin{eqnarray}\label{eq_mu3}
\dys \mu(B_r(x)) & = & \lim_{\eps\to0}\mu_\eps(B_r(x)) \nonumber \\
\nonumber\\
& = & \dys \lim_{\eps\to0} F_\eps(\ue, \Om\cap B_r(x), \partial\Om\cap B_r(x)) \nonumber \\
\nonumber\\
& \geq & \dys\liminf_{\eps\to0} (1-\delta(\Psi_r))^{p+3}F_\eps({\bar{u}_{\eps}}, D_r, E_r).
\end{eqnarray}
Moreover, by Proposition \ref{pro_slice}, we have
\begin{equation}\label{eq_mu3.2}
\dys \liminf_{\eps\to0} F_\eps({\bar{u}_{\eps}}, D_r, E_r)\geq\gamma_p\left|\int_{S\bar{v}\cap E_r}\nu_v d\hs^1\right|.
\end{equation}
Finally, we notice that $\delta(\Psi_r)$ vanishes and $\dys \left|\int_{S\bar{v}\cap E_r}\nu_v d\hs^1\right|=2r+o(r)$ as $r$ goes to 0. So (\ref{eq_mu3}) and (\ref{eq_mu3.2}) give the following inequality
$$
\dys \frac{\mu(B_r(x))}{2r}\geq\gamma_p\left(1+\frac{o(r)}{2r}\right) \ \text{as} \ r\to0,
$$
that implies $\mu\geq\mu^3$. This concludes the proof of the lower bound inequality. \hfill $\Box$

\bigskip

\subsection{Upper bound inequality}

We will construct an optimal sequence $(\ue)$ according to Theorem \ref{c4_teoremap}-(iii) in a suitable partition of $\Om$, as in \cite[Theorem 2.6-(iii)]{alberti98}, but for the estimate of the boundary effect we will use the optimal profile problem (\ref{profilep2D}) in connection with the results proved in the previous section. We do not present all the details, just sketch the main ideas and state the needed lemmas.

Fix $(u,v) \in BV(\Om;\{\al',\be'\})\times BV(\partial\Om; \{\al',\be'\})$. It is not restrictive to assume that the singular sets $Su$ and $Sv$ are closed manifolds of class $C^2$ without boundary (see \cite[Theorem 1.24]{giusti}). We may also assume that $u$ and $v$ (up to modifications on negligible sets) are constant in each connected component of $\Omega\setminus Su$ and $\partial\Omega\setminus Sv$, respectively.

\medskip

The idea is to construct a {\em partition} of $\Omega$ in four subsets, and to use the preliminary convergence results of the previous sections to obtain the upper bound inequality.
\begin{figure}[htp!]
\centering
\psfrag{ua}{${u\!=\!\al}$} \psfrag{ub}{${u\!=\!\be}$} \psfrag{vba}{${v\!=\!\al'}$} \psfrag{vbp}{$v\!=\!\be'$} \psfrag{a1}{\!\!$A_1$} \psfrag{a2}{\!\!$A_2$} \psfrag{b1}{\!\!$B_1$} \psfrag{b2}{\!\!$B_2$} \psfrag{su}{\!\!\!\!$Su\cap\Gamma_r$} \psfrag{sv}{\!$Sv$}  
\includegraphics[scale=0.34]{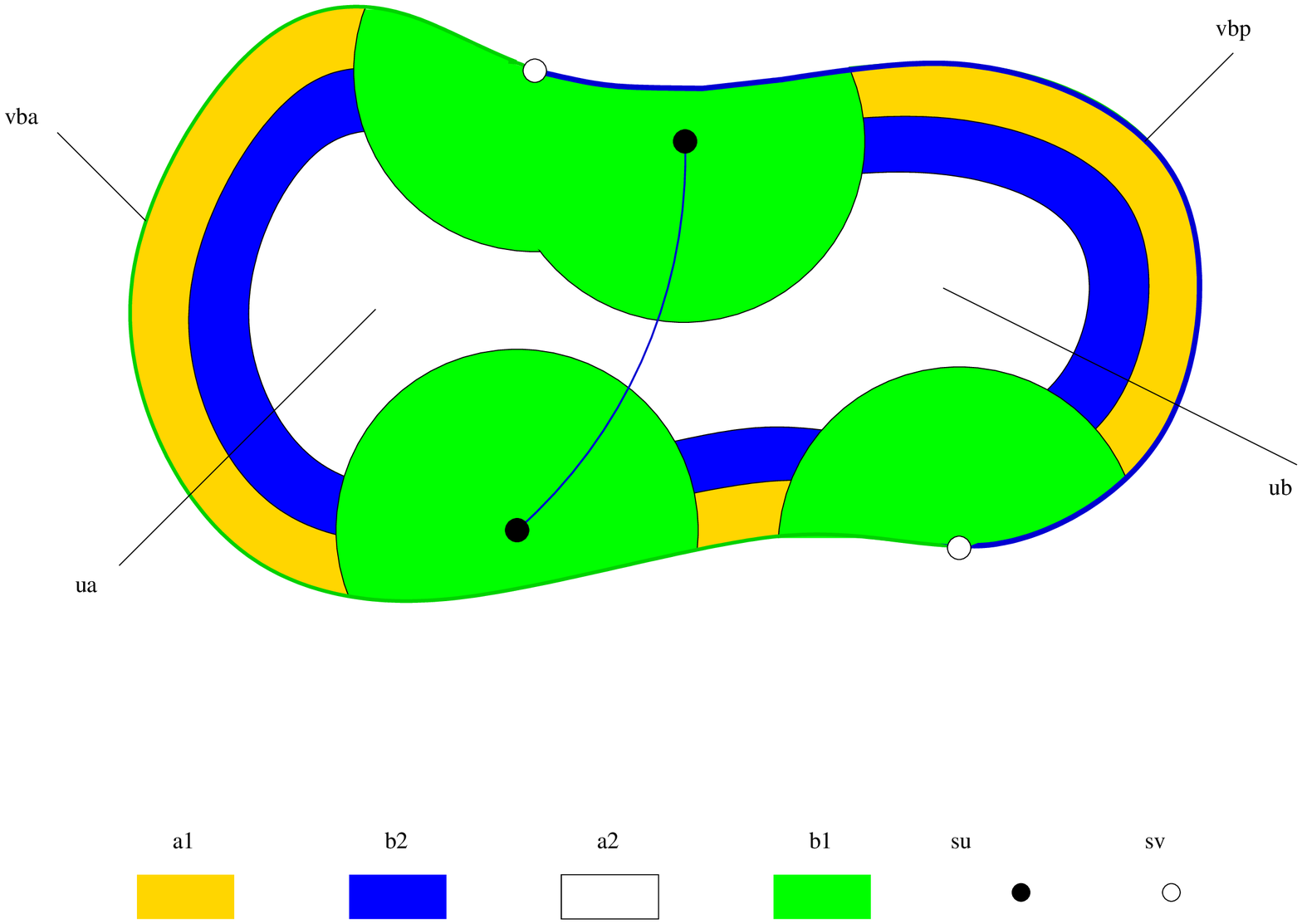}
\caption{Upper bound inequality - partition of $\Om$ (see \cite[Fig. 6]{alberti98}).}\label{fig_partition}
\end{figure}
\noindent
\\ For every $r>0$, we set
\begin{equation*}
\dys \Gamma_r:=\left\{x\in \Omega : \text{dist}(x,\partial\Omega)=r\right\}.
\end{equation*}

\noindent
\\ {\it Step 1~: Construction of the partition}.
\noindent
\\ Fix $r>0$ such that $\Gamma_r$ and $\Gamma_{2r}$ are Lipschitz surfaces and $Su\cap\Gamma_r$ is a Lipschitz curve. 

Now, we are ready to construct the following partition of $\Omega$:
\begin{eqnarray*}
B_1 \!\!&:=&\!\!\left\{x\in \Omega : \text{dist}(x, Sv\cup(Su\cap\Gamma_r))<3r\right\}, \\
 A_1\!\!&:=&\!\!\left\{x\in\Omega\setminus \overline{B}_{1}: \text{dist}(x,\partial\Omega)<r\right\}, \\
 B_2\!\!&:=&\!\!\left\{x\in\Omega\setminus \overline{B}_{1}: r<\text{dist}(x,\partial\Omega)<2r\right\}, \\
 A_2\!\!&:=&\!\!\left\{x\in\Omega\setminus \overline{B}_{1}: \text{dist}(x,\partial\Omega)>2r\right\}.
\end{eqnarray*}
(See Fig. \ref{fig_partition})

\medskip

For every $r>0$ and every $\dys \eps<r^{\frac{p-1}{p-2}}$ 
 we construct a Lipschitz function $\ue=u_{\eps,r}$ in each subset, with controlled Lipschitz constant.

\bigskip

\noindent
\\ {\it Step 2~: Construction of $u_{\eps,r}$ in $A_2$}.
\noindent
\\ In $A_2$, we take $\ue$ 
being the optimal sequence for the functional $G_\eps$ in the set $A_2$ as in Theorem \ref{teo_modicap}-(iii) 
and we extend it to $\partial A_2$ by continuity.
Note that $\ue$ converges to $u$ pointwise in $A_2$ and uniformly on $\partial A_2\cap \partial B_2$, and
\begin{eqnarray}\label{ub_a2}
\dys F_\eps(\ue, A_2, \emptyset)  \equiv G_\eps(\ue, A_2) \! & \leq\! & \sigma_p\hs^2(Su\cap A_2)+ o(1) \nonumber \\
\\
& \leq\! & \sigma_p\hs^2(Su\cap A_2)+ o(1) , \ \ \text{as} \ \eps \to 0. \nonumber
\end{eqnarray}

\bigskip

\noindent
\\ {\it Step 3~: Construction of $u_{\eps,r}$ in $A_1$}.
\noindent
\\ The function $u$ is constant (equal to $\alpha$ or $\beta$) on every connected component $A$ of $A_1$, and the function $v$ is constant (equal to $\alpha'$ or $\beta'$) on $\partial A\cap\partial\Omega$.
We can extend it to $\partial A_1$ with continuity; Proposition \ref{pro_modica2p}-(ii) 
gives 
\begin{eqnarray}\label{ub_a1}
\dys F_\eps(\ue, A_1, \partial A_1\cap\partial\Omega) \! \equiv \! G_\eps(\ue, A_1) \, & \leq & \, c_p\int_{\partial A_1\cap\partial\Omega}|\mathcal{W}(Tu(x))-\mathcal{W}(v(x))|d\hs^{2} \nonumber \\
\\
& & + o(1), \ \ \text{as} \ \eps \to 0. \nonumber
\end{eqnarray}

\bigskip

\noindent
\\ {\it Step 4~: Construction of $u_{\eps,r}$ in $B_2$}.
\noindent
\\ Following \cite{alberti98}, to construct $\ue$ on $B_2$, we need to ``glue'' the values of $A_1$ and $A_2$. Take a cut-off function $\xi$ such that $\xi=1$ in $A_1$ and $\xi=0$ in $A_2$ and consider the function
$$
\ue=\xi\bar{u}_1+(1-\xi)\bar{u}_2,
$$
where $\bar{u}_i$ is the extension to $B_2$ of $\ue|_{A_i}$. Then, when $\eps\to 0$, by the decay of the function $\xi$, we have
\begin{eqnarray*}
\dys \eps^{p-2}\int_{B_2}|D\ue|^p h^{2-p}dx \leq C\left(\int_{B_2}|D\bar{u}_1|^p h^{2-p}dx+\int_{B_2}|D\bar{u}_2|^p h^{2-p}dx \right. \hskip 2cm\\
\left. +\int_{B_2}|D\xi|^p|\bar{u}_1-\bar{u}_2|^p h^{2-p}dx\right) = o(1).
\end{eqnarray*}

\bigskip

\noindent
\\ {\it Step 5~: Construction of $u_{\eps,r}$ in $B_2$}.
\noindent
\\ Finally, for the last part $B_1$, we will use an optimal profile for the minimum problem (\ref{profilep2D}).

By Proposition \ref{pro_min2d}, there exists $\psi \in L^1_{\text{loc}}(\R^2_+)$ such that $T\psi(t)\to\alpha'$ as $t\to -\infty$, $T\psi(t)\to\beta'$ as $t\to +\infty$ and $\dys H_1(\psi, \R^2_+, \R)=\gamma_p.$
We can construct a function $\we:\R^2_+\to\R$ following the method used to provide a good competitor $u_\delta$ in the proof of Proposition \ref{pro_lb_disco}.

For every $\eps>0$, $\roe, \sie \in \R$, we take a cut-off function $\xi\in C^{\infty}(\R^2_+)$ such that $\xi\equiv 1$ on ${(\R^2_+)\setminus D_{\roe}}$ and $\xi\equiv 0$ on ${D_{\sie}}$ such that $|{D}\xi|\leq \frac{1}{|\roe-\sie|} $. We denote by $\bar{u}$ the function defined in polar coordinates $\theta \in [0,\pi)$, $\rho \in [0,+\infty)$ as follows
$$
\dys \bar{u}(\theta, \rho):=\frac{\theta}{\pi}\alpha'+\left(1-\frac{\theta}{\pi}\right)\beta'.
$$
We define $\we$ as
$$
\we(x):=
\begin{cases}
\psi(\frac{x}{\eps}) & \text{if} \ x\in D_{\sie}, \\\\
\xi(x)\bar{u}(x)+(1-\xi(x))\psi(\frac{x}{\eps}) & \text{if} \ x \in D_{\roe}\setminus D_{\sie}, \\\\
\bar{u}(x) & \text{if} \ x \in \R^2_{+} \setminus D_{\roe}.
\end{cases}
$$
Let us show that we can choose $\roe$ and $\sie$ such that $\we$ satisfies the following inequality
\begin{equation}\label{ub_heps}
\eps^{p-2}\int_{D_{\roe}}|D\we|^p h^{2-p}dx +\frac{1}{\sqrt{\eps}}\int_{E_{\roe}}V(T\we) \leq \gamma_p+o(1), \ \text{as} \ \eps\to 0.
\end{equation}

\smallskip

 By the definition of $\we$ and by standard changing variable formula ($y=x/\sqrt{\eps}$), we have
\begin{eqnarray}\label{eq_b1}
\dys \eps^{p-2}\int_{D_{\rho_{\eps}}}|D\we|^p h^{2-p} dx \! & = & \! \int_{(D_{\rho_{\eps}}\cap D_1^0)/\sqrt{\eps}}|D\we^{(\eps)}|^p y_2^{2-p} dy \nonumber \\
& & \! + \eps^{p-2}\int_{D_{\rho_{\eps}}\cap(D_1^0)^{c}}|D\we|^p \left(1-\sqrt{x_1^2+x^2_2}\right)^{2-p} dx  \nonumber \\
\nonumber \\
& \leq & \! \int_{\R_+^2}|D\psi|^p y_2^{2-p} dy + \eps^{p-2}\int_{(D_{\rho_{\eps}}\setminus D_{\sie})\cap D_1^0}|D\we|^p x_2^{2-p} dx  \\
& & \! + \eps^{p-2}\int_{D_{\rho_{\eps}}\cap(D_1^0)^{c}}|D\we|^p \left(1-\sqrt{x_1^2+x^2_2}\right)^{2-p} dx  \nonumber \\
& =: & \! \int_{\R_+^2}|D\psi|^p y_2^{2-p} dy + I_1 + I_2, \nonumber
\end{eqnarray}
where, $D^1_0$ is defined by (\ref{def_d10}). \smallskip 

Notice that when $\rho_\eps <<1$, the integral $I_2$ is zero, since $D_{\roe} \cap {(D_{1}^{0})}^c$ is empty, hence we have
\begin{eqnarray}\label{eq_2eto}
\dys H_\eps(\we, D_{\roe}, E_{\roe}) \! & \leq & \! H_1(\psi, \R^2_+, \R) + I_1 \nonumber \\
\\
& = & \gamma_p +I_1. \nonumber
\end{eqnarray}
Thus, to obtain (\ref{ub_heps}), it suffices to estimate the integral $I_1$. We may work more or less like in the proof of Proposition \ref{pro_lb_disco}.\smallskip

We have \begin{eqnarray}\label{eq_2etoetoeto}
\dys I_1 \! & \leq & \! 3^{p-1} \eps^{p-2}\int_{(D_{\roe}\setminus D_{\sie})\cap D_1^0}|D\psi(\frac{x}{\eps})|^p x_2^{2-p} dx + 3^{p-1} \eps^{p-2}\int_{(D_{\roe}\setminus D_{\sie})\cap D_1^0}|D\bar{u}|^p x_2^{2-p} dx \nonumber \\
\\
& & \! + 3^{p-1} \eps^{p-2}\int_{(D_{\roe}\setminus D_{\sie})\cap D_1^0}|D\xi|^p|\psi(\frac{x}{\eps})-\bar{u}|^p x_2^{2-p} dx \nonumber
\end{eqnarray}
and the last two integrals in the right part of (\ref{eq_2etoetoeto}) can be explicitly estimated as follows
\begin{eqnarray}\label{eq_3eto}
\dys 3^{p-1}\eps^{p-2}\int_{(D_{\roe}\setminus D_{\sie})\cap D_1^0}|D\bar{u}|^p x_2^{2-p} dx \! & = & \! 3^{p-1}\eps^{p-2}\frac{|\be'-\al'|^p}{\pi^p}\int_0^{\pi}\!\!\int_{\sie}^{\roe}\frac{(\rho\sin{\theta})^{2-p}}{\rho^p} \rho d\rho d\theta\nonumber \\
\\
& \leq & \! C_1\frac{\eps^{p-2}}{{\roe}^{2(p-2)}} \nonumber
\end{eqnarray}
and
\begin{eqnarray}\label{eq_3etoeto}
\dys 3^{p-1}\eps^{p-2}\!\!\int_{D_{\roe}\setminus D_{\sie}\cap (D_1^0)^c}\!\!|D\xi|^p|\psi(\frac{x}{\eps})-\bar{u}|^p x_2^{2-p} dx \! \! & \leq & \!\! 3^{p-1}\eps^{p-2}\!\frac{(2m)^p}{(\roe-\sie)^p}\int_0^\pi \!\!\int_{\sie}^{\roe}\!\!\rho(\rho\sin{\theta})^{2-p} d\theta d\rho \nonumber \\
\\
& \leq & \! \! C_2\frac{\eps^{p-2}\roe^{4-p}}{(\roe-\sie)^p}. \nonumber
\end{eqnarray}
\smallskip

Finally, by (\ref{eq_2etoetoeto}), (\ref{eq_3eto}) and (\ref{eq_3etoeto}), the inequality (\ref{eq_2eto}) becomes
\begin{eqnarray}\label{eq_4eto}
\dys H_\eps(\we, D_{\roe}, E_{\roe}) \! & \leq & \! \gamma_p + 3^{p-1} \eps^{p-2}\int_{(D_{\roe}\setminus D_{\sie})\cap D_1^0}|D\psi(\frac{x}{\eps})|^p x_2^{2-p} dx \nonumber \\
\nonumber \\
& & \! + C_1\frac{\eps^{p-2}}{\roe^{2(p-2)}} + C_2\frac{\eps^{p-2}\roe^{4-p}}{(\roe-\sie)^p} \nonumber \\
\\
& \leq & \! \gamma_p + o(1) \ \ \ \text{as} \ \eps \to 0, \nonumber
\end{eqnarray}
where we also used that, since $\dys \int_{\R_+^2}\!|D\psi|^p y_2^{2-p}$ is finite, by suitable choosing $\roe$ and $\sie$ we get
\begin{eqnarray}\label{eq_4bis}
3^{p-1} \eps^{p-2}\int_{(D_{\roe}\setminus D_{\sie})\cap D_1^0}|D\psi(\frac{x}{\eps})|^p x_2^{2-p} dx  \! & = & \! 3^{p-1}\int _{((D_{\roe}\setminus D_{\sie})\cap D_1^0)/\sqrt{\eps}}|D\psi|^p y_2^{2-p} dy \nonumber \\
\\
& = & \! o(1) \ \ \ \text{as} \ \eps \to 0. \nonumber
\end{eqnarray}
\smallskip

Since the neighborhood $B_1$ is Lipschitz equivalent (modulo some multiplicative constant) to the product $Sv \times D_{\roe}$, we can construct the following transplanted function $\bar{w}_\eps$
$$
\dys \bar{w}_\eps(x,z):=\we(x) \ \ \ \forall x\in Sv, \forall z\in \R^2_+.
$$
By Fubini's Theorem, we obtain 
\begin{eqnarray}\label{eq_5eto}
F_\eps(\bar{w}_\eps, Sv\!\times\! D_{\roe}, Sv\!\times\! E_{\roe}) \!\!\!  & = & \!\! \!\hs^1(Sv)\!\left(\!H_\eps(\we, D_{\roe}, E_{\roe})+\frac{1}{\eps^{\frac{p-2}{p-1}}}\!\!\int_{D_{\roe}}\!\!W(\we)h^{\frac{p-2}{p-1}}dx \!\right) \nonumber\\
\\
& \leq & \!\!\! \hs^1(Sv)\left(\!H_\eps(\we, D_{\roe}, E_{\roe})+C_3\frac{\roe^2}{\eps^{\frac{p-2}{p-1}}}\right). \nonumber 
\end{eqnarray} \smallskip
Hence, by suitably choosing $\roe$ and $\sie$(i.e., such that $\frac{\eps^{p-2}}{\roe^{2(p-2)}}$, $\frac{\eps^{p-2}\roe^{4-p}}{(\roe-\sie)^p}$ and $\frac{\roe^2}{\eps^{\frac{p-2}{p-1}}} \to 0$ as $\eps \to 0$ and (\ref{eq_4bis}) holds; for instance, $\dys \roe=\eps^{\frac{p-2}{p-1}}$ and $\sie=\eps^{\frac{p-2}{2(p-1)}}$),
we get
\begin{equation}\label{eq_40star}
\dys F_\eps(\bar{w}_\eps, Sv\!\times\! D_{\roe}, Sv\!\times\! E_{\roe})\leq \hs^1(Sv)\left(\gamma_p+o(1) \right) \ \ \text{as} \ \eps \to 0.
\end{equation}
\medskip

\noindent
\\ {\it Step 6~: The upper bound inequality}.
\noindent
\\ Now, we can use an extension lemma for the remaining pieces, which is contained in \cite[Lemma 5.4]{palatucci07}.

\begin{lemma}\label{lem_lip} Let $A$ be a domain in $\R^3, A'\subset \partial A, v:A'\to[-m, m]$ a Lipschitz function {\rm(}where $m$ is given by {\rm (\ref{ipotesi}))} and $G_\eps$ defined by {\rm (\ref{fu_momow})}.
\\ Then, for every $\eps > 0$, there exists an extension $u:\overline{A}\to[-m, m]$ such that
$$
\text{\rm Lip}(u)\leq{\eps^{-\frac{p-2}{p-1}}}+\text{\rm Lip}(v)
$$
and
\begin{equation}\label{eq_lip}
G_\eps(u, A)\leq\left( (\eps^{\frac{p-2}{p-1}}\text{\rm Lip}(v)+1)^p+C_m\right)\left(\hs^2(\partial A)+o(1)\right)\omega, \ \ \text{as} \ \eps \to 0,
\end{equation}
where $\dys C_m:=\!\!\max_{t\in [-m, m]}\!\! W(t)$, $\omega:=\min\{\|v-\alpha\|_{L^\infty},\|v-\beta\|_{L^\infty}\}.$
\end{lemma}
\medskip
\noindent

The rest of the proof of the theorem follows one of the author \cite{palatucci07} or Alberti, Bouchitt\'e and Seppecher \cite{alberti98} with minor modifications, and we can find a Lipschitz function $\ue$ in the whole $\Om$ with the required behavior.  \hfill $\Box$

 \bigskip

\end{document}